\documentclass[twoside, 10pt]{article}

\usepackage{amsmath}
\usepackage{paralist}
\usepackage{amssymb,mathrsfs,psfrag,graphicx}
\usepackage{hyperref}
\usepackage{CJK,color}

\allowdisplaybreaks

\definecolor{red}{rgb}{1.00,0.00,0.00}
\definecolor{blue}{rgb}{0.00,0.00,0.63}
\definecolor{black}{rgb}{0.00,0.00,0.00}

\newcommand{\nm}{\nonumber}

\newcommand\tbbint{{-\mkern -16mu\int}}

\newcommand\dbbint{{-\mkern -19mu\int}}

\newcommand\bbint{
{\mathchoice{\dbbint}{\tbbint}{\tbbint}{\tbbint}}
}

\def\v{\varepsilon}

\def\t{\theta}
\def\T{\Theta}
\def\k{\kappa}

\def\mb{\mathbf}
\def\a{\alpha}

\def\d{\delta}

\def\f{\frac}

\def\z{\zeta}

\def\di{\displaystyle}

\def\pa{\partial}

\newtheorem{theorem}{Theorem}[section]

\newtheorem{lemma}[theorem]{Lemma}
\newtheorem{proposition}{Proposition}[section]

\newtheorem{remark}{Remark}

\newcommand{\ba}{\begin{aligned}}
\newcommand{\ea}{\end{aligned}}
\newcommand{\be}{\begin{equation}}
\newcommand{\ee}{\end{equation}}

\begin{document}

\title{\bf   Vanishing dissipation limit for non-isentropic Navier-Stokes equations with shock data} \vskip 0.5cm
\author{
 \quad Feimin Huang\thanks{Academy of Mathematics and Systems Science, Chinese Academy of Sciences, Beijing 100190, P. R. China and School of Mathematical Sciences, University of Chinese Academy of Sciences, Beijing 100049, P. R. China
   ({\tt fhuang@amt.ac.cn}). },
\quad Teng Wang\thanks{Corresponding author. School of Mathematics, Statistics and Mechanics, Beijing University of Technology, Beijing 100124, P. R. China
   ({\tt tengwang@amss.ac.cn}). }
 }

\date{}
\maketitle
\begin{abstract}

This paper is concerned with the vanishing dissipation limiting problem of one-dimensional non-isentropic Navier-Stokes equations with shock data. The limiting problem was solved in 1989 by Hoff-Liu in \cite{HL-1} for isentropic gas with single shock, but was left open for non-isentropic case. In this paper, we solve the non-isentropic case, i.e., we first establish the global existence of
solutions to the non-isentropic Navier-Stokes equations with initial discontinuous shock data, and then show these solutions converge in $L^{\infty}$ norm to a single shock wave of the corresponding Euler equations away from the shock curve in any finite time interval, as both the viscosity and heat-conductivity tend to zero.

Different from \cite{HL-1} in which an integrated system was essentially  used, motivated by \cite{KV-JEMS, KV-Invention}, we introduce a time-dependent shift $\mb{X}^\varepsilon(t)$  to the viscous shock so that a weighted Poincar\'{e} inequality can be applied to overcome the difficulty generated from the ``bad" sign of the derivative of viscous shock velocity, and the anti-derivative technique is not needed. We also obtain an intrinsic property of non-isentropic viscous shock, see Lemma \ref{le-im} below. With the help of Lemma \ref{le-im}, we can derive the desired uniform a  priori estimates of solutions, which can be shown to converge in $L^{\infty}$ norm to a single inviscid shock in any given finite time interval away from the shock, as the vanishing dissipation limit. Moreover,  the shift $\mb{X}^\varepsilon(t)$ tends to zero in any finite time as viscosity tends to zero.  The proof consists of a scaling argument, $L^2$-contraction technique with time-dependent shift to the shock, and relative entropy method.

\


\end{abstract}



%
%
\section{Introduction}
\renewcommand{\theequation}{\arabic{section}.\arabic{equation}}
\setcounter{equation}{0}

In this paper, we are concerned with the one-dimensional non-isentropic compressible Navier-Stokes system in Lagrangian coordinates:
\begin{equation}\label{ns}
\left\{
\ba
& v_t-u_x=0,\\
& u_t+p_x=\varepsilon\left(\frac{u_x}{v}\right)_x,\\
& \left(e+\frac{u^2}{2}\right)_t+\left(p u\right)_x=\left(\kappa \frac{\t_x}{v}+\v\frac{uu_x}{v}\right)_x
\ea
\right.
\end{equation}
for $t>0$, $x\in\mathbb{R}=(-\infty,+\infty)$, where the unknown functions $v=v(t,x)>0$, $u=u(t,x)$ and $\t=\t(t,x)>0$ represent the specific volume, fluid velocity and an absolute temperature, respectively. Here we consider the polytropic gas, that is, the pressure $p$ and the internal energy $e$ are explicitly given as functions of the specific volume $v$ and the absolute temperature $\theta$:
$$
p=\frac{R\t}{v}=Av^{-\gamma}\exp\left(\f{\gamma-1}{R}S\right),\quad e=\f{R\t}{\gamma-1}+\mathrm{constant},
$$
where $S$ is the entropy, $\gamma>1$ is the adiabatic exponent, and $A$ and $R$ are both positive constants. In addition, $\v>0$ is the viscosity constant and $\k>0$ is the coefficient of heat-conductivity. We consider the following Riemann initial data
\begin{equation}\label{R-in}
(v,u,\t)(0,x)=\left\{
\ba
& (v_-,u_-,\t_-),&  x<0,\\
& (v_+,u_+,\t_+),&  x>0,
\ea
\right.
\end{equation}
where $v_{\pm}(>0)$, $u_{\pm}$ and $\t_{\pm}(>0)$ are given constants.

As in \cite{JNS}, we assume the
viscosity $\v$ and the heat-conductive coefficient $\k$ satisfy
\be\label{relation}
\left\{
\ba
& \k=O(\v), & {\rm as} \quad \v\rightarrow0; \\
& \f{\k(\v)}{\v}\geq c>0, &{\rm as}\quad \v\rightarrow0,
\ea
\right.
\ee
for some positive constant $c$. Then the non-isentropic compressible Navier-Stokes system \eqref{ns} formally converges to the corresponding non-isentropic compressible Euler system as $\v\rightarrow 0+$,
\begin{equation}\label{euler}
\left\{
\ba
& v_t-u_x=0,\\
& u_t+p_x=0,\\
& \left(e+\frac{u^2}{2}\right)_t+\left(p u\right)_x=0.
\ea
\right.
\end{equation}

A fundamental and challenging problem in mathematics is whether  the vanishing dissipation limit of the solution of the Navier-Stokes equations to the Euler equations could be proved  rigorously. Because the inviscid Euler system \eqref{euler} usually contains discontinuities, such as shock waves and contact discontinuities, it is difficult to show the limit in the general setting, especially for the discontinuous initial data. Essential new ideas and methods are needed to handle this problem. Therefore, any attempt on this problem that involves the singularity in the inviscid solution can be viewed as progress towards this general program.

There is an extensive literature on the vanishing viscosity limit to the inviscid compressible flow in one space dimension. For a system of hyperbolic conservation laws with artificial viscosity, Goodman-Xin \cite{Goodman-Xin} studied the limit to the piecewise smooth solutions having a finite number of noninteracting shocks. Later, Yu \cite{Yu-ARMA} improved the results of \cite{Goodman-Xin} to allow initial layers by a detailed pointwise analysis. As a breakthrough result, Bianchini-Bressan \cite{BB}
proved the vanishing viscosity limit to a strictly hyperbolic $n\times n$ system with small BV initial data. It is noted that the work of \cite{BB} essentially depends on the viscosity matrix  $\v\mathbb{I}$. There are other interesting works on the
inviscid limit for small BV solution, see \cite{B,BY,BHWY} and the references therein.  For $L^{\infty}$ entropy solution, which is less regular than the BV solution,
the vanishing artificial viscosity limit to the isentropic Euler equations with large $L^{\infty}$ initial data has been established by DiPerna \cite{Diperna}, Chen \cite{Chen-acta,Chen-proc,CP}, Lions-Perthame-Souganidis \cite{Lions-PS}, Lions-Perthame-Tadmor \cite{Lions-PT}, Huang-Wang \cite{Huang-WangZ},  Schrecker-Schulz \cite{SS} via the method of compensated compactness. See also Ding-Chen-Luo \cite{DCL} for numerical approximation limit. However, little information on the $L^{\infty}$ entropy solutions is known through the compensated compactness.

The study on the vanishing physical viscosity was initiated by Stokes \cite{stokes}. The first rigorous limit of viscous shock layer from the Navier-Stokes system to the Euler system was obtained by Gilbarg \cite{Gilbarg}. Hoff-Liu \cite{HL-1} in 1989 first proved the vanishing viscosity limit for piecewise constant shock.
Indeed, they showed that the solution of the isentropic Navier-Stokes system with initial shock data exists all the time and converges to the inviscid single shock as the viscosity vanishes, uniformly away from the shock.  Zhang-Pan-Tan \cite{ZPT} extended the result of \cite{HL-1} to a composite wave consisting of two shocks. However, the limiting problem toward a single shock was left open for non-isentropic case until now.

In the present paper, we solve the limiting problem mentioned above, that is, the vanishing dissipation limit to a single shock for non-isentropic Navier-Stokes system \eqref{ns} with initial shock data \eqref{R-in}. More precisely, we prove that the solution of the system \eqref{ns} exists for all time with discontinuous initial data \eqref{R-in}, and converges in $L^{\infty}$ norm to a single inviscid shock on any given finite time interval  away from the shock, as the viscosity and heat-conductive coefficients tend to zero. This gives a rigorous mathematical justification of the limit for the non-isentropic compressible Navier-Stokes system in the presence of single shock wave.


We now introduce some preliminary notations and give some background materials before stating the main theorem. The inviscid system \eqref{euler} has three distinct eigenvalues for positive $v$ and $\t$
$$
\lambda_1(v,\t)=-\f{\sqrt{\gamma R\t}}{v},\quad \lambda_2=0, \quad\lambda_3(v,\t)=\f{\sqrt{\gamma R\t}}{v},
$$
which implies the first and third characteristic fields are genuinely nonlinear and
the second field is linearly degenerate.  In this paper, we focus our attention on the single shock case, that is, the entropy solution of inviscid system \eqref{euler} with \eqref{R-in} is a single shock wave.

We firstly recall the shock wave curve in the phase plane $(v,u,\theta)$.
Given the right end state $(v_+,u_+,\t_+)$, the $i-$shock wave curve of the Euler equations \eqref{euler} for the left end state $(v,u,\t)$ is defined as follows:

$\bullet$ $i-$shock wave curve $(i=1,3)$,
$$
\ba
S_i(v_+,u_+,\t_+):=\left\{(v,u,\t)   \middle |
\ba
&-s_i(v_+-v)-(u_+-u)=0,\\
&-s_i(u_+-u)+(p_+-p)=0,\quad {\rm and}\quad\lambda_{i}^+<s_i<\lambda_{i}^-,\\
&-s_i\left[e_++\f{u_+^2}{2}-\left(e+\f{u^2}{2}\right)\right]+(p_+u_+-pu)=0,
\ea
\right\},
\ea
$$
where $e_+=\f{R\t_+}{\gamma-1} +$constant, $p_+=\f{R\t_+}{v_+}$, $\lambda_{i}^{\pm}=\lambda_i(v_\pm,\t_\pm)$, and $s_i$ is the $i-$shock speed.

It is known that  the solution to the Riemann problem \eqref{euler}, \eqref{R-in} admits a $3-$shock wave if $(v_-,u_-,\t_-)\in S_3(v_+,u_+,\t_+)$ (cf. \cite{Lax}) and the explicit formula of $3-$shock wave is given as
\be\label{discon-shock}
(v^S,u^S,\t^S)(t,x):=\left\{
\ba
& (v_-,u_-,\t_-), & x<s_3 t, \\
& (v_+,u_+,\t_+), & x>s_3 t,
\ea
\right.
\ee
where $(v_-,u_-,\t_-)$, $(v_+,u_+,\t_+)$ and the shock speed $s_3$ are constants satisfying the Rankine-Hugoniot (R-H) condition
\be\label{RH}
\left\{
\ba
&-s_3(v_+-v_-)-(u_+-u_-)=0,\\
&-s_3(u_+-u_-)+(p_+-p_-)=0,\\
&-s_3\left[e_++\f{u_+^2}{2}-\left(e_-+\f{(u_-)^2}{2}\right)\right]+(p_+u_+-p_-u_-)=0,
\ea
\right.
\ee
and Lax' shock condition that
\be\label{Lax-condi}
\lambda_3^+<s_3<\lambda_3^-
\ee
with $\lambda_3^+=\lambda_3(v_+,\t_+)$ and $\lambda_3^-=\lambda_3(v_-,\t_-)$.
Note that \eqref{Lax-condi} is equivalent to
$$
u_+<u_-\quad ({\rm or}\,\,  s_3(v_+-v_-)>0).
$$
We denote the strength of $3-$shock wave by
$
\d=|p_--p_+|,
$
where $p_\pm=\f{R\t_\pm}{v_\pm}$. It follows from the R-H condition \eqref{RH} that $\d$ is equivalent to $|v_--v_+|+|u_--u_+|+|\t_--\t_+|$.

We are ready to state our main theorem on the vanishing dissipation limit of solution for the non-isentropic compressible Navier-Stokes equations \eqref{ns} with initial shock data \eqref{R-in}.

\begin{theorem}\label{thm}
For any given $(v_+,u_+,\t_+)$, suppose $(v_-,u_-,\t_-)\in S_3(v_+,u_+,\t_+)$. Let $( v^S, u^S, \t^S)(t,x)$ be the $3-$shock wave defined in \eqref{discon-shock} and the  viscosity and heat-conductivity satisfy \eqref{relation}. Then there exist positive constants $\d_0$, $\varepsilon_0$ independent of time $t$, such that if $\v\leq \v_0$ and the wave strength $\d\leq \d_0$, the initial value problem \eqref{ns}, \eqref{R-in} admits a unique global piecewise smooth solution $(v^{\v},u^{\v},\t^{\v})(t,x)$, for each $\v>0$, satisfying the following properties:

$\bullet$ $u^{\v}$ and $\t^{\v}$ are continuous for $t>0$; $u^{\v}_x$, $\t^{\v}_x$, $v^{\v}$, $v^{\v}_x$ and $v^{\v}_t$ are uniformly H\"{o}lder continuous in the sets $\{(t,x)\,|\, t\geq \tau ,x<0\}$ and $\{(t,x)\,|\,  t\geq \tau,x>0\}$ for any $\tau>0$; $u^{\v}_t$, $u^{\v}_{xx}$, $v^{\v}_{xt}$, $\t^{\v}_t$ and $\t^{\v}_{xx}$ are H\"{o}lder continuous on compact set in $\{(t,x)\,|\, t>0, x\neq0\}$. Moreover, the jumps in $v^{\v}(t,x)$, $u^{\v}_x(t,x)$ and $\t^{\v}_x(t,x)$ at $x=0$ satisfy
$$
|[v^{\v}(t,0)]|, \quad |[u_x^{\v}(t,0)]|,\quad |[\t_x^{\v}(t,0)]|\leq Ce^{-\f{ct}{\v}},
$$
where $C$ and $c$ are positive constants independent of $t$ and $\v$.

$\bullet$  Let $T$ be a positive constant, then the solutions $(v^{\v},u^{\v},\t^{\v})(t,x)$ converge to the $3-$shock wave $(v^S,u^S,\t^S)(t,x)$ up to time $T$, as  $\v\rightarrow0$, uniformly away from $t=0$ and the shock, i.e., for any positive number $h$, we have
\be\label{limit}
\lim_{\v\rightarrow0}\sup_{h\leq t \leq T}\sup_{|x-s_3t|\geq h}|(v^{\v},u^{\v},\t^{\v})(t,x)-(v^S,u^S,\t^S)(t,x)|=0.
\ee
\end{theorem}


\begin{remark}
 Hoff-Liu \cite{HL-1} in 1989 first investigated the inviscid limit to a single shock for isentropic Navier-Stokes equations. However, the corresponding result for non-isentropic case remains open until now. Our main result in Theorem \ref{thm} solves
the limiting problem to a single shock for non-isentropic case.
\end{remark}

\begin{remark}
Our result still holds for the case $1-$shock, since the proof for $1-$shock is identical.
\end{remark}

Before sketching our strategy to prove Theorem \ref{thm}, we give a brief recall of the interesting work of Hoff-Liu \cite{HL-1}. Since the anti-derivative technique is used, the initial perturbation in the anti-derivative form  $\|\int_{-\infty}^x(U-\Phi)(y,0)dy\|_{L^2}$  is around $\d^{-1/2}$, which is large as $\delta$ is small, and this brings new challenge in the small perturbation framework. Here $U$ is the solution of Navier-Stokes equations and $\Phi$ is the traveling wave. A key observation in \cite{HL-1} is that  $\|\int_{-\infty}^x(U-\Phi)(y,t)dy\|_{L^2}^2$  is of the order $\d^{a_1}(t+1)^{b_1}$ (where $a_1-2b_1\geq0$) and $\|U-\Phi\|_{L^2}^2$ is of the order $\d^{a_2}(t+1)^{b_2}$ (where $a_2-2b_2\geq1$), and this leads to the desired smallness of the $L^{\infty}$-norm of $\int_{-\infty}^x(U-\Phi)(y,t)dy$ at $t=\d^{-2-\vartheta}$ so that the small perturbation framework is still available on $t>\d^{-2-\vartheta}$.


However, the above method in \cite{HL-1} can not be applied directly to the non-isentropic flow \eqref{ns} since  $\int_{-\infty}^x(U-\Phi)(t,y)dy$ may not belong to $L^2$ space. In fact, the vanishing dissipation limit problem to a single shock for the non-isentropic Navier-Stokes system has been open for long time. In this paper,
we solve this open problem.
Motivated by the $L^2$ contraction method developed in \cite{KV-JEMS,KV-Invention},
we study the limiting problem in the original $H^1$-perturbation framework instead of the anti-derivative technique.
Different from \cite{KV-JEMS,KV-Invention,KVW} for the isentropic system, it is hard to introduce the new effective velocity
for the non-isentropic system \eqref{ns} since the pressure $p$ is a function of both the specific volume $v$ and the absolute temperature $\t$. Instead, we  decompose the ``bad" term in \eqref{s3-eta} in form of the velocity $u$ by Taylor expansion  with small viscosity and wave strength (see \eqref{s3-eta}). Then, we derive a uniform estimate in the new variable on the variable transformation. The key point is to establish Lemma \ref{le-im} based on the underlying structure of non-isentropic viscous shock wave, which is quite different from the isentropic case, where the pressure is a function only in the specific volume $v$. Thus, with the help of the  inequality \eqref{important}, we can obtain the $L^2$-contraction around viscous shock profile for non-isentropic case by introducing a suitable shift $\mb{X}^{\v}(t)$ and a weighted  Poincar\'{e} inequality \eqref{poin-2}. Moreover, the shift $\mb{X}^{\v}(t)$ is shown to tend to zero in any finite time interval as the viscosity tends to zero. Thus we can
prove that the solutions of the non-isentropic Navier-Stokes equations \eqref{ns} converge to a single inviscid shock as the viscosity and heat-conductive coefficients tend to zero, away from the shock.

There are other interesting works on the inviscid limit of compressible Navier-Stokes equations in the presence of basic wave patterns, see \cite{HJW,HWY-KRM,HWY-ARMA,HWWY,JNS,KV-Invention,Ma,WY,Xin,Xin-Z,ZPWT}, and the references cited therein. In particular, Kang-Vasseur \cite{KV-Invention} proved the stability and uniqueness of inviscid shock in the relative entropy form for $p$-system in the weak sense through the vanishing viscosity limit from the Navier-Stokes system.

Recently, the $L^2$-contraction method is also employed in \cite{KVW} to prove the time asymptotic stability of composite waves with a viscous shock and a rarefaction wave to the isentropic Navier-Stokes equations, which solves an important open problem proposed by Matsumura-Nishihara \cite{MN-92} in 1992. We refer to \cite{KVW-2} for the case of a viscous shock, a viscous contact wave, and an inviscid rarefaction wave to the non-isentropic Navier-Stokes system \eqref{ns}.

The rest part of the paper is organized as follows. In Sect. \ref{shock}, we construct the smooth viscous shock wave and show its properties. In Sect. \ref{pf-thm}, we prove our main vanishing dissipation limit result by continuation arguments. In Sect. \ref{pf-priori}, we prove a-priori estimates in Proposition \ref{priori}.

{\textbf{Notation.}} In the paper, we always use the notation $\bbint_{\mathbb{R}}=\int_{\mathbb{R}_-}+\int_{\mathbb{R}_+}$,  $\|\cdot\|$ to denote the usual $L^2(\mathbb{R})$ norm, $\| \cdot \rlap{--}\|$ to denote the piecewise $L^2(\mathbb{R}_{\pm})$ norm, that is, $\|f\rlap{--}\|^2=\bbint_{\mathbb{R}}f^2dx$. $\|\cdot\|_{H^1}$ and $\|\cdot\rlap{--}\|_{H^1}$ represent the $H^1(\mathbb{R})$ and $H^1(\mathbb{R}_{\pm})$ norm, respectively. And the notation $[f]$ represents the jump of the function $f$ at $x=0$ or $y=0$ if without confusion.


\section{Viscous shock wave}\label{shock}
\setcounter{equation}{0}

Now we recall the definitions of viscous shock waves of \eqref{ns} which correspond to the shock wave \eqref{discon-shock}. We see that the $3-$viscous shock wave which corresponds to the $3-$shock wave is a traveling wave solution of \eqref{ns} with the formula $(V^{\v},U^{\v},\T^{\v})(\xi)$ which is determined by
\be\label{shock-equ}
\left\{
\ba
&-s_3(V^{\v})'-(U^{\v})'=0,\\
&-s_3(U^{\v})'+(P^{\v})'=\v\left(\f{(U^{\v})'}{V^{\v}}\right)',\qquad\,\, ':=\f{d}{d\xi},\qquad\xi=x-s_3t,\\
&-s_3\left(\f{R\T^{\v}}{\gamma-1}+\f{(U^{\v})^2}{2}\right)'+(P^{\v}U^{\v})'=\k\left(\f{(\T^{\v})'}{V^{\v}}\right)'+\v\left(\f{(U^{\v})' U^{\v}}{V^{\v}}\right)',
\ea
\right.
\ee
and
\be\label{shock-far}
(V^{\v},U^{\v},\T^{\v})(-\infty)=(v_-,u_-,\t_-),\qquad (V^{\v},U^{\v},\T^{\v})(+\infty)=(v_+,u_+,\t_+),
\ee
where $P^{\v}=P(V^{\v},\T^{\v})=\f{R\T^{\v}}{V^{\v}}$, and $(v_\pm,u_\pm,\t_{\pm})$ satisfy the R-H condition \eqref{RH} and Lax entropy condition \eqref{Lax-condi}, and the shock speed $s_3$ is uniquely determined by \eqref{RH}. The existence and uniqueness (up to a shift in $\xi$) of the solution of \eqref{shock-equ}, \eqref{shock-far} can be found in \cite{KM}.


The property of the $3-$viscous shock wave $(V^{\v}, U^{\v},\T^{\v})(x-s_3t)$ are summarized in the following lemma, which can be found in \cite{Huang-matsumura,KM}.

\begin{lemma}[\cite{Huang-matsumura,KM}]\label{le-shock}
Let $(V^{\v},U^{\v},\T^{\v})$ be the smooth $3-$viscous shock wave with the strength $\d=|p_--p_+|$
\begin{itemize}
\item [(i)]
$U^{\v}_x=-s_3V^{\v}_x<0$,\quad $s_3\T^{\v}_x<0$,\quad $x\in\mathbb{R}$,\quad $t\geq0$,

\item [(ii)]  There exist positive constants $c$ and $C$ which depend only on $(v_-,u_-,\t_-)$ and $(v_+,u_+,\t_+)$ such that
$$
\ba
& |(V^{\v}-v_-,U^{\v}-u_-,\T^{\v}-\t_-)|\leq C\d e^{-c\d\f{|x-s_3t|}{\v}}, & x<s_3t, \quad t\geq 0,\\
& |(V^{\v}-v_+,U^{\v}-u_+,\T^{\v}-\t_+)|\leq C\d e^{-c\d\f{|x-s_3t|}{\v}}, & x>s_3t, \quad t\geq 0,\\
& |(V^{\v}_x,U^{\v}_x,\T^{\v}_x)|\leq C\f{\d^2}{\v}e^{-c\d\f{|x-s_3t|}{\v}}, & x\in\mathbb{R}, \quad
t\geq 0,\\
& |(V^{\v}_{xx},U^{\v}_{xx},\T^{\v}_{xx})|\leq C\f{\d}{\v}|(V^{\v}_{x},U^{\v}_{x},\T^{\v}_{x})|, & x\in\mathbb{R}, \quad
t\geq 0.
\ea
$$
\item [(iii)] There exist positive constants $c$ and $C$ which depend only on $(v_-,u_-,\t_-)$ and $(v_+,u_+,\t_+)$ such that
$$
\Big |s_3^2-\f{\gamma p_-}{v_-}\Big|\leq C\d,\qquad \Big |s_3^2-\f{\gamma p_+}{v_+}\Big|\leq C\d.
$$
\end{itemize}

\end{lemma}

\textbf{\emph{Proof}}: The estimates in (i) and (ii) can be shown by direct calculations, we omit the proof. In the following, we prove (iii). Integrating \eqref{shock-equ} over $\mathbb{R}$, we have
$$
\ba
0&=-\f{R}{\gamma-1}(\t_+-\t_-)-p_-(v_+-v_-)+\f{s_3^2}{2}(v_+-v_-)^2\\
&=-\f{R\t_-}{\gamma-1}\Big(\f{\t_+}{\t_-}-1\Big)-R\t_-\Big(\f{v_+}{v_-}-1\Big)
+\f{s_3^2v_-^2}{2} \Big(\f{v_+}{v_-}-1\Big)^2,
\ea
$$
which leads to
\be\label{add-2}
\f{\t_+}{\t_-}-1=-(\gamma-1)\Big(\f{v_+}{v_-}-1\Big)+\f{\gamma-1}{R\t_-}\f{s_3^2v_-^2}{2}\Big(\f{v_+}{v_-}-1\Big)^2.
\ee
Note that the R-H condition \eqref{RH} and \eqref{add-2} give
$$
\ba
s_3^2&=-\f{p_+-p_-}{v_+-v_-}=-\f{p_-}{v_+}\Big(\f{\t_+}{\t_-}-\f{v_+}{v_-}\Big) \Big/ \Big(\f{v_+}{v_-}-1\Big)
=\f{p_-}{v_+}\Big[1-\Big(\f{\t_+}{\t_-}-1\Big) \Big/ \Big(\f{v_+}{v_-}-1\Big)\Big]\\
&=\f{p_-}{v_+}\Big[1-\Big(-(\gamma-1)+\f{\gamma-1}{R\t_-}\f{s_3^2v_-^2}{2}\Big(\f{v_+}{v_-}-1\Big)\Big)\Big]
=\f{\gamma p_-}{v_+}-\f{\gamma-1}{2v_+}s_3^2(v_+-v_-),
\ea
$$
which yields (iii).

%

\hfill $\Box$

\begin{remark}
In \cite{Huang-matsumura,KM}, it requires $1<\gamma\leq 2$. In fact, as pointed in \cite{Huang-matsumura}, this condition might be removed by more careful and complicated estimates.
\end{remark}

The following lemma plays a crucial role to obtain the extra weak dissipation effect along the shock propagation direction, which fully used the compressible structure of the non-isentropic  viscous shock wave.

\begin{lemma}\label{le-im}
Let $(V^{\v},U^{\v},\T^{\v})$ be the $3-$viscous shock wave with strength $\d=|p_--p_+|$, then there exists positive constant $C$ such that
\be\label{important}
\left|\f{V^{\v}-v_-}{P^{\v}-p_-}-\f{V^{\v}-v_+}{P^{\v}-p_+}-\f{\gamma+1}{2\gamma p_+}\f{R\v\gamma}{\k(\gamma-1)^2+R\v\gamma}(v_+-v_-)\right|\leq C\d^2,
\ee
where $p_\pm=\f{R\t_\pm}{v_\pm}$ and $P^{\v}=P(V^{\v},\T^{\v})=\f{R\T^{\v}}{V^{\v}}$.
\end{lemma}

\textbf{\emph{Proof}}: Integrating \eqref{shock-equ} over  $[\xi,+\infty)$ or $(-\infty,\xi]$, we have
\be\label{integ-t}
\ba
\f{\k\T^{\v}_x}{s_3V}
&=-\f{R}{\gamma-1}(\T^{\v}-\t_\pm)-p_\pm(V^{\v}-v_\pm)+\f{s_3^2}{2}(V^{\v}-v_\pm)^2\\
&=-\f{R\t_\pm}{\gamma-1}\Big(\f{\T^{\v}}{\t_\pm}-1\Big)-R\t_\pm\Big(\f{V^{\v}}{v_\pm}-1\Big)
+\f{s_3^2v_\pm^2}{2}\Big(\f{V^{\v}}{v_\pm}-1\Big)^2,
\ea
\ee
which together with (iii) in Lemma \ref{le-shock} gives
\be\label{T+-}
\ba
\f{\T^{\v}}{\t_\pm}-1&=-(\gamma-1)\Big(\f{V^{\v}}{v_\pm}-1\Big)+\f{\gamma-1}{R\t_\pm}
\f{s_3^2v_\pm^2}{2}\Big(\f{V^{\v}}{v_\pm}-1\Big)^2-\f{\gamma-1}{R\t_\pm}\f{\k\T^{\v}_x}{s_3V^{\v}}\\
&=-(\gamma-1)\Big(\f{V^{\v}}{v_\pm}-1\Big)+\f{\gamma(\gamma-1)}{2}\Big(\f{V^{\v}}{v_\pm}-1\Big)^2
-\f{\gamma-1}{R\t_\pm}\f{\k\T^{\v}_x}{s_3V^{\v}}+O(\d^3).
\ea
\ee
Then, using \eqref{T+-}, we have
\be\label{pv}
\ba
\f{V^{\v}-v_-}{P^{\v}-p_-}-\f{V^{\v}-v_+}{P^{\v}-p_+}
&=\f{V^{\v}}{p_-} \Big(\f{V^{\v}}{v_-}-1\Big) \Big/ \Big( \f{\T^{\v}}{\t_-}-\f{V^{\v}}{v_-}\Big)
-\f{V^{\v}}{p_+} \Big(\f{V^{\v}}{v_+}-1\Big) \Big/ \Big( \f{\T^{\v}}{\t_+}-\f{V^{\v}}{v_+}\Big)\\
&=\f{V^{\v}}{p_-}\f{1}{\Big(\f{\T^{\v}}{\t_-}-1\Big) \Big/  \Big(\f{V^{\v}}{v_-}-1\Big)-1}
-\f{V^{\v}}{p_+}\f{1}{\Big(\f{\T^{\v}}{\t_+}-1\Big) \Big/  \Big(\f{V^{\v}}{v_+}-1\Big)-1}\\
&=\Big(\f{V^{\v}}{p_-}-\f{V^{\v}}{p_+}\Big)\f{1}{\Big(\f{\T^{\v}}{\t_-}-1\Big) \Big/  \Big(\f{V^{\v}}{v_-}-1\Big)-1}\\
&\quad+\f{V^{\v}}{p_+}\left[\f{1}{\Big(\f{\T^{\v}}{\t_-}-1\Big) \Big/  \Big(\f{V^{\v}}{v_-}-1\Big)-1}-\f{1}{\Big(\f{\T^{\v}}{\t_+}-1\Big) \Big/  \Big(\f{V^{\v}}{v_+}-1\Big)-1}\right]\\
&=:I_1+I_2.
\ea
\ee
It follows from \eqref{add-2} that
\be\label{p-p+}
\ba
p_- - p_+&=R\Big(\f{\t_-}{v_-}-\f{\t_+}{v_+}\Big)
=\f{R\t_-}{v_+}\Big(\f{v_+}{v_-}-\f{\t_+}{\t_-}\Big)
=\f{R\t_-}{v_+}\Big[\Big(\f{v_+}{v_-}-1\Big)
-\Big(\f{\t_+}{\t_-}-1\Big)
\Big]\\
&=\f{R\t_-}{v_+}\gamma\Big(\f{v_+}{v_-}-1\Big)+O(\d^2)
=\f{\gamma p_+}{v_+}(v_+-v_-)+O(\d^2).
\ea
\ee
Then, using \eqref{T+-} and \eqref{p-p+}, we have
\be\label{I1}
\ba
I_1&=\f{V^{\v}(p_+-p_-)}{p_-p_+}\Big/
\Big(-\gamma+\f{\gamma(\gamma-1)}{2}\Big(\f{V^{\v}}{v_-}-1\Big)
-\f{\gamma-1}{p_-}\f{\k\T^{\v}_x}{s_3 V^{\v}}\f{1}{V^{\v}-v_-}+O(\d^2)\Big)\\
&=\f{v_+}{\gamma p_+^2}(p_- -p_+)+O(\d^2)
=\f{1}{p_+}(v_+-v_-)+O(\d^2).
\ea
\ee
It follows from \eqref{T+-} and \eqref{p-p+} that
$$
\ba
I_2&=\f{V^{\v}}{p_+}\Big[1 \Big/ \Big(-\gamma+\f{\gamma(\gamma-1)}{2}\Big(\f{V^{\v}}{v_-}-1\Big)
-\f{\gamma-1}{p_-}\f{\k\T^{\v}_x}{s_3 V^{\v}}\f{1}{V^{\v}-v_-}+O(\d^2)
\Big)\\
&\qquad-1 \Big/ \Big(-\gamma+\f{\gamma(\gamma-1)}{2}\Big(\f{V^{\v}}{v_+}-1\Big)
-\f{\gamma-1}{p_+}\f{\k\T^{\v}_x}{s_3 V^{\v}}\f{1}{V^{\v}-v_+}+O(\d^2)
\Big)
\Big]\\
&=\f{V^{\v}}{\gamma^2p_+}\Big[\f{\gamma(\gamma-1)}{2}V^{\v}\Big(\f{1}{v_+}-\f{1}{v_-}\Big)
+\f{\gamma-1}{p_+}\f{\k\T^{\v}_x}{s_3 V^{\v}}\Big(\f{1}{V^{\v}-v_-}-\f{1}{V^{\v}-v_+}\Big)\\
&\quad+(\gamma-1)\Big(\f{1}{p_-}-\f{1}{p_+}\Big)\f{\k\T^{\v}_x}{s_3 V^{\v}}\f{1}{V^{\v}-v_-}
\Big]+O(\d^2)\\
&=-\f{\gamma-1}{2\gamma p_+}(v_+-v_-)+\f{(\gamma-1)v_+}{\gamma^2p_+^2}
\f{\k\T^{\v}_x}{s_3 V^{\v}}\Big(\f{1}{V^{\v}-v_-}-\f{1}{V^{\v}-v_+}\Big)+O(\d^2),
\ea
$$
using \eqref{integ-t}, \eqref{p-p+} and (iii) in Lemma \ref{le-shock}, we have
$$
\ba
&\f{\k\T^{\v}_x}{s_3 V^{\v}}\Big(\f{1}{V^{\v}-v_-}-\f{1}{V^{\v}-v_+}\Big)\\
=&\f{R}{\gamma-1}\Big(\f{\T^{\v}-\t_+}{V^{\v}-v_+}-\f{\T^{\v}-\t_-}{V^{\v}-v_-}\Big)
+(p_+-p_-)+\f{s_3^2}{2}(v_+-v_-)\\
=&\f{R\v\gamma^2-2\k\gamma(\gamma-1)}{\k(\gamma-1)^2+R\v\gamma}\f{p_-}{2v_-}(v_+-v_-)
-\f{\gamma p_+}{v_+}(v_+-v_-)+\f{\gamma p_+}{2v_+}(v_+-v_-)+O(\d^2)\\
=&\f{R\v\gamma^2-2\k\gamma(\gamma-1)}{\k(\gamma-1)^2+R\v\gamma}\f{p_+}{2v_+}(v_+-v_-)
-\f{\gamma p_+}{2v_+}(v_+-v_-)+O(\d^2),
\ea
$$
where we have used the fact (see formula (B.5) of \cite{KVW-2}) that
$$
\ba
\f{\T^{\v}-\t_+}{V^{\v}-v_+}-\f{\T^{\v}-\t_-}{V^{\v}-v_-}=\f{\gamma-1}{2R}
\f{R\v\gamma^2-2\k\gamma(\gamma-1)}{\k(\gamma-1)^2+R\v\gamma}\f{p_-}{v_-}(v_+-v_-)+O(\d^2).
\ea
$$
Thus we have
\be\label{I2}
I_2=-\f{\gamma-1}{\gamma p_+}(v_+-v_-)
+\f{\gamma-1}{2\gamma p_+}\f{R\v\gamma-2\k(\gamma-1)}{\k(\gamma-1)^2+R\v\gamma}(v_+-v_-)+O(\d^2).
\ee
Finally, we substitute \eqref{I1} and \eqref{I2} into \eqref{pv} to obtain the desired inequality \eqref{important}, and the proof of Lemma \ref{le-im} is completed.

\hfill $\Box$

Finally, we introduce a one-dimensional weighted sharp Poincar\'{e} type inequality, whose proof can be found in \cite{KV-JEMS} and plays a very important role in our stability analysis.
\begin{lemma}[Lemma 2.9 \cite{KV-JEMS}]
For any $f:[0,1]\rightarrow\mathbb{R}$ satisfying $\int_0^1z(1-z)|f'|^2dz<\infty$, it holds
\be\label{poin-2}
\int_0^1\left|f-\int_0^1fdz\right|^2dz\leq \f12\int_0^1z(1-z)|f'|^2dz.
\ee
\end{lemma}

%

%
%
%
%


%
%
%
%

\section{Proof of the main result}\label{pf-thm}
\setcounter{equation}{0}

In this section, we will prove our main result, Theorem \ref{thm}. Introduce the following scaled variables
\be\label{scaling}
y=\f{x}{\v},\quad \tau=\f{t}{\v},
\ee
and set
$$
(v^{\v},u^{\v},\t^{\v})(t,x)=(v,u,\t)(\tau,y).
$$
Then the new unknown function $(v,u,\t)(\tau,y)$ satisfies the system
\begin{equation}\label{new-euler}
\left\{
\ba
& v_\tau-u_y=0,\\
& u_\tau+p_y=\left(\f{u_y}{v}\right)_y,\\
& \f{R}{\gamma-1}\t_{\tau}+p u_y=\nu\left(\f{\t_y}{v}\right)_y+\f{u_y^2}{v}
\ea
\right.
\end{equation}
with initial data
\be\label{new-data}
(v,u,\t)(0,y)=\left\{
\ba
& (v_-,u_-,\t_-),&  y<0,\\
& (v_+,u_+,\t_+),&  y>0,
\ea
\right.
\ee
where $\nu=\f{\k}{\v}$ denotes the coefficient of scaled heat conductivity and satisfies $\underline{\nu}\leq\nu\leq \bar\nu$, uniformly in $\v\rightarrow 0+$, for some positive constants $\underline{\nu}$ and $\bar \nu$.

\subsection{Construction of the shift}

Here and in what follows, we introduce the notation: for any function $f:\mathbb{R}^+\times\mathbb{R}\rightarrow\mathbb{R}$ and the shift $\mb{X}(\tau)$,
$$
f^{-\mb{X}}:=f(\tau,y-\mb{X}(\tau)).
$$
The shift $\mathbf{X}(\tau)$ is defined as a solution to the ODE:
\be\label{X}
\left\{
\ba
\dot{\mathbf{X}}(\tau)&=-\f{M}{\d}\int_{\mathbb{R}}
a^{-\mathbf{X}}\left(
-\f{P^{-\mb{X}} V^{-\mathbf{X}}_y}{s_3 V^{-\mb{X}}}+U^{-\mathbf{X}}_y
+\f{p_+\T_y^{-\mathbf{X}}}{s_3 \T^{-\mb{X}}}
\right)(u-U^{-\mb{X}})dy,\\
\mathbf{X}(0)&=0
\ea
\right.
\ee
with the constant $M=\f{(\gamma+1)^2}{2\a_+v_+^2}\f{2\nu(\gamma-1)^2+R\gamma}{\nu(\gamma-1)^2+R\gamma}$. Here $a(\cdot):\mathbb{R}\rightarrow\mathbb{R}$ is a weight function defined by
\be\label{a}
a(y-s_3\tau)=1+\f{p_--P(y-s_3\tau)}{\sqrt{\d}}=1+\f{p(v_-,\t_-)-P(V,\T)(y-s_3\tau)}{\sqrt{\d}}.
\ee
Obviously, the weight function $a$ satisfies
\be\label{a>1}
1\leq a\leq 1+\sqrt{\d},\qquad  a_y=-\f{P_y}{\sqrt{\d}}>0.
\ee





Let $F(\tau,\mathbf{X}(\tau))$ be the right-hand side of the ODE $\eqref{X}_1$. By Lemma \ref{le-shock} and \eqref{a>1}, we know that $a$, $a_y$, $U$ are bounded, and $V$, $\T$ are bounded from below and above. Indeed, from the Propositions \ref{local}, \ref{priori}, we can also show $u$ is bounded. Thanks to these facts,  there is some constant $C>0$ such that
\be\label{F(X)}
\sup_{\mathbf{X}\in\mathbb{R}}|F(\tau,\mathbf{X})|\leq\f{C}{\d}\|a\|_{L^{\infty}}
(\|u\|_{L^{\infty}}+\|U^{-\mb{X}}\|_{L^{\infty}})
\left\|\left(V_y^{-\mb{X}},U_y^{-\mb{X}},\T_y^{-\mb{X}}\right)\right\|_{L^1}\leq C,
\ee
and
\be
\ba
\sup_{\mathbf{X}\in\mathbb{R}}|\pa_{\mathbf{X}}F(\tau,\mathbf{X})|&\leq \f{C}{\d}\|a_y\|_{L^{\infty}}(\|u\|_{L^{\infty}}+\|U^{-\mb{X}}\|_{L^{\infty}})
\left\|\left(V_y^{-\mb{X}},U_y^{-\mb{X}},\T_y^{-\mb{X}}\right)\right\|_{L^1}\\
&+\f{C}{\d}\|a\|_{L^{\infty}}(\|u\|_{L^{\infty}}+\|U^{-\mb{X}}\|_{L^{\infty}})
\left\|\left(V_{yy}^{-\mb{X}},U_{yy}^{-\mb{X}},\T_{yy}^{-\mb{X}}\right)\right\|_{L^1}\\
&+\f{C}{\d}\|a\|_{L^{\infty}}
\left\|\left(V_y^{-\mb{X}},U_y^{-\mb{X}},\T_y^{-\mb{X}}\right)U_y^{-\mb{X}}\right\|_{L^1}
\leq C.
\ea
\ee
Then ODE \eqref{X} has a unique absolutely continuous solution $\mathbf{X}(\tau)$ defined on any interval in time $[0,\mathcal{T}]$ by the well-known Cauchy-Lipschitz theorem.
In particular, since $|\dot{\mathbf{X}}(\tau)|\leq C$ by \eqref{F(X)}, we can obtain
$$
|\mathbf{X}(\tau)|\leq C\tau, \qquad \forall \tau\in[0, \mathcal{T}].
$$


\subsection{Reformulation of the problem}

Setting
$$
\ba(\phi,\psi,\z)(\tau,y)&=(v,u,\t)(\tau,y)-(V^{-\mathbf{X}},U^{-\mathbf{X}},\T^{-\mathbf{X}})(\tau,y)\\
&=(v,u,\t)(\tau,y)-(V,U,\T)(y-s_3\tau-\mathbf{X}(\tau)),
\ea
$$
where $(V,U,\T)$ satisfies \eqref{shock-equ} with $\v=1$ and $\k=\nu$. Then, it follows from \eqref{new-euler}, \eqref{new-data} and \eqref{shock-equ} with $\v=1$ and $\k=\nu$ that
\be\label{perturb}
\left\{
\ba
&\phi_\tau-\dot{\mathbf{X}}(\tau)V^{-\mathbf{X}}_y-\psi_y=0,\\
&\psi_\tau-\dot{\mathbf{X}}(\tau)U^{-\mathbf{X}}_y+(p-P^{-\mathbf{X}})_y
=\left(\f{u_y}{v}-\f{U^{\mathbf{-X}}_y}{V^{-\mathbf{X}}}\right)_y,\\
&\f{R\z_{\tau}}{\gamma-1}-\dot{\mathbf{X}}(\tau)\f{R\T^{-\mathbf{X}}_y}{\gamma-1}+p\psi_y
+(p-P^{-\mathbf{X}})U^{-\mathbf{X}}_y\\
&=\nu\left(\f{\t_y}{v}-\f{\T^{-\mathbf{X}}_y}{ V^{-\mathbf{X}}}\right)_y+\f{\psi_y^2}{v}+2\f{\psi_y U^{-\mathbf{X}}_y}{v}+ (U^{-\mathbf{X}}_y)^2\left(\f{1}{v}-\f{1}{ V^{-\mathbf{X}}}\right),
\ea
\right.
\ee
and the initial value
\begin{equation}\label{per-initial}
\ba
(\phi,\psi,\z)(0,y)=(\phi_0,\psi_0,\z_0)(y)=\left\{
\ba
& (\,v_- -V(y),u_- - U(y),\t_- -\T(y)\,), &\di y<0,\\
& (\,v_+ -V(y),u_+ - U(y),\t_+ -\T(y)\,), &\di y>0.
\ea
\right.
\ea
\end{equation}
We can see that the initial data $(\phi_0,\psi_0,\z_0)(y)$ and its derivatives are sufficiently smooth away from but up to $y=0$ and
$$
(\phi_0,\psi_0,\z_0)(y)\in L^2(\mathbb{R}),\quad \phi_{0y}\in L^2(\mathbb{R}_{\pm}).
$$
For simplicity, we denote
$$
\mathcal{N}_0:=\|(\phi_0,\psi_0,\z_0)\|+\|\phi_{0y}\rlap{--}\|.
$$
We show that the initial value problem \eqref{perturb}, \eqref{per-initial} has a unique global solution $(\phi,\psi,\z)(\tau,y)$, with regularity as in Theorem \ref{thm}, which approaches zero uniformly as $\tau\rightarrow\infty$.

\begin{proposition}\label{pro-1}
There exists a positive constant $\d_0$ such that if the initial data and the wave strength $\d$ satisfies
$$
\mathcal{N}_0+\d\leq \d_0,
$$
then the problem \eqref{perturb}, \eqref{per-initial} admits a unique global solution $(\phi,\psi,\z)(\tau,y)$ satisfying

\begin{itemize}
\item [(i)] There exists a positive constant $C$ independent of $\tau$ such that
    \be\label{pro-es1}
    \sup_{\tau\geq0}\left(\|(\phi,\psi,\z)(\tau,\cdot)\|^2
    +\|\phi_y(\tau,\cdot)\rlap{--}\hspace{0.07em}\|^2\right)    +\int_0^{\infty}\left(\|(\psi_y,\z_y)(\tau,\cdot)\|^2
    +\|\phi_y(\tau,\cdot)\rlap{--}\hspace{0.07em}\|^2\right)d\tau
    \leq C(\mathcal{N}_0^2+\d).
    \ee

\item[(ii)] For any $\tau_*>0$, there exists a positive constant $C=C(\tau_*)$ depending on $\tau_*$, such that
    \be\label{pro-es2}
    \sup_{\tau\geq\tau_*}\|(\psi_y,\z_y,\psi_\tau,\z_\tau)(\tau,\cdot)\rlap{--}\hspace{0.07em}\|^2
    +\int_{\tau_*}^{\infty}\|(\psi_{yy},\z_{yy},\psi_{\tau y},\z_{\tau y})(\tau,\cdot)\rlap{--}\hspace{0.07em}\|^2d\tau
    \leq C(\tau_*)(\mathcal{N}_0^2+\d).
    \ee

\item[(iii)] The jump condition of $\phi(\tau,y)$ at $y=0$ admits the bound
\be\label{pro-jump}
|[\phi](\tau)|\leq C \d e^{-c\tau},
\ee
where the positive constants $C$ and $c$ are independent of $\tau\in(0,\infty)$.


\end{itemize}
\end{proposition}

Assume the Proposition \ref{pro-1} holds, then for any $\tau_*>0$, we have from the inequalities \eqref{pro-es1}, \eqref{pro-es2}, and equations \eqref{perturb} that
$$
\int_{\tau_*}^{\infty}\left(\|(\phi_y,\psi_y,\z_y)\rlap{--}\|^2
+\left|\f{d}{d\tau}\|(\phi_y,\psi_y,\z_y)\rlap{--}\|^2\right|\right)d\tau<\infty,
$$
which leads to
$$
\lim_{\tau\rightarrow\infty}\|(\phi_y,\psi_y,\z_y)\rlap{--}\|^2=0,
$$
from which and Gagliardo-Nirenberg inequality, it holds
$$
\ba
\lim_{\tau\rightarrow\infty}\sup_{y\neq0}\|(\phi,\psi,\z)(\tau)\|_{L^{\infty}}
&\leq C\lim_{\tau\rightarrow\infty}\|(\phi,\psi,\z)(\tau)\|^{\f12}
\|(\phi_y,\psi_y,\z_y)(\tau)\rlap{--}\|^{\f12}\\
&\leq C\lim_{\tau\rightarrow\infty}\|(\phi_y,\psi_y,\z_y)(\tau)\rlap{--}\|^{\f12}=0.
\ea
$$
This together with the estimate \eqref{pro-jump} gives
\be\label{limit-1}
\lim_{\tau\rightarrow\infty}\sup_{y\in\mathbb{R}}
\left|\left(v-V^{-\mb{X}},u-U^{-\mb{X}},\t-\T^{-\mb{X}}\right)(\tau,y)\right|=0.
\ee

Next, it follows from \eqref{X} and Lemma \ref{le-shock} that
$$
|\dot{\mathbf{X}}(\tau)|\leq \f{C}{\d}\|\psi\|_{L^{\infty}}\|(V^{-\mathbf{X}}_y,U^{-\mathbf{X}}_y,\T^{-\mathbf{X}}_y)\|_{L^1}
\leq C\|\psi\|_{L^{\infty}},
$$
which with $\mathbf{X}(0)=0$ and H\"{o}lder inequality yield
$$
\ba
|\mb{X}(\tau)|&\leq C\int_0^{\tau}\|\psi\|_{L^{\infty}}d\tau\leq C\int_0^{\tau}\|\psi\|^{\f12}\|\psi_y\|^{\f12}d\tau
\leq C\left(\int_0^{\tau}\|\psi_y\|^2d\tau\right)^{\f14}
\left(\int_0^{\tau}\|\psi\|^{\f23}d\tau\right)^{\f34}\\
&\leq C\left(\int_0^{\tau}\|\psi_y\|^2d\tau\right)^{\f14}\sup_{\tau>0}\|\psi\|^{\f12}\tau^{\f34}.
\ea
$$
Recalling the scaled variable \eqref{scaling}, the above inequality indicates
$$
\left|\v\mathbf{X}\left(\f{t}{\v}\right)\right|\leq C\varepsilon^{\f14}T^{\f34},
$$
for $t\in[0,T]$ with any finite positive number $T$. Then, for any positive number $h$, there exists sufficient small $\v$ such that
\be\label{es-X-2}
\left|\v\mathbf{X}\left(\f{t}{\v}\right)\right|\leq \f{h}{2}.
\ee
%
Therefore, we have
\be\label{limit-2}
\ba
&(v^{\v},u^{\v},\t^{\v})(t,x)-(v^S,u^S,\t^{S})(t,x)\\
=&(v,u,\t)(\tau,y)-(v^S,u^S,\t^{S})(\tau,y)\\
=&(v,u,\t)(\tau,y)
-(V^{-\mb{X}},U^{-\mb{X}},\T^{-\mb{X}})(\tau,y)
+(V^{-\mb{X}},U^{-\mb{X}},\T^{-\mb{X}})(\tau,y)
-(v^S,u^S,\t^{S})(\tau,y),
\ea
\ee
where the second last term in the last line of \eqref{limit-2} vanishes as $\tau\rightarrow+\infty$ due to the large time behavior \eqref{limit-1}, which means it vanishes as  $\v\rightarrow0+$,   uniformly with $t\geq h$ for any positive number $h$ due to a straightforward scaling argument $\tau=\f{t}{\v}$.

Then for the last term in the last line of \eqref{limit-2}, using Lemma \ref{le-shock} and \eqref{es-X-2}, for any $h>0$, it holds
\be\label{limit-3}
\ba
&\lim_{\v\rightarrow 0}\sup_{|y-s_3\tau|\geq \f{h}{\v}}|(V^{-\mb{X}},U^{-\mb{X}},\T^{-\mb{X}})(\tau,y)
-(v^S,u^S,\t^{S})(\tau,y)|\\
=&\lim_{\v\rightarrow 0}\sup_{|y-s_3\tau|\geq \f{h}{\v}}|(V,U,\T)(y-s_3\tau-\mb{X}(\tau))
-(v^S,u^S,\t^{S})(\tau,y)|\\
=&\lim_{\v\rightarrow 0}\sup_{|x-s_3t|\geq h}\left|(V^{\v},U^{\v},\T^{\v})\left(\f{x-s_3t}{\v}-\mb{X}\Big(\f{t}{\v}\Big)\right)
-(v^S,u^S,\t^{S})(t,x)\right|\\
\leq& \lim_{\v\rightarrow 0}\sup_{|x-s_3t|\geq h} C\d \exp
\left(-\f{c\d}{\v}\left|x-s_3t-\v\mb{X}\Big(\f{t}{\v}\Big) \right| \right)
\leq \lim_{\v\rightarrow 0}C\d \exp\left(-\f{c\d h}{2\v}\right)=0.
\ea
\ee
Thus combination of \eqref{limit-2} and \eqref{limit-3} gives the dissipation limit behavior \eqref{limit} of the solutions $(v^{\v},u^{\v},\t^{\v})$. Then the proof of Theorem \ref{thm} can be finished with the regularity results from Proposition \ref{pro-1}.

\hfill $\Box$

Denote

$$
\ba
N(\tau_*,\tau^*)&=\sup_{\tau\in[\tau_*,\tau^*]}\Big(\|(\phi,\psi,\z)(\tau,\cdot)\|
+\|(\phi_y,\psi_y,\z_y)(\tau,\cdot)\rlap{--}\|+\|(\psi_\tau,\z_{\tau})(\tau,\cdot)\|\Big),\\
N(\tau_*)&=N(\tau_*,\tau^*),
\ea
$$
and define the solution space of the problem \eqref{perturb},  \eqref{per-initial}, as follows
$$
\ba
\mathcal{X}[\tau_*,\tau^*]=\Big\{(\phi,\psi,\z)\,| \,
&(\phi,\psi,\z)\in C([\tau_*,\tau^*];H^1(\mathbb{R}_{\pm})), \,
\phi_y\in L^2(\tau_*,\tau^*;L^2(\mathbb{R}_{\pm})), \\
&  (\psi_y,\z_y)\in L^2(\tau_*,\tau^*;H^1(\mathbb{R}_{\pm})),\\
&(\psi_\tau,\z_\tau)\in L^{\infty}(\tau_*,\tau^*;L^2(\mathbb{R}_{\pm}))\cap
L^2(\tau_*,\tau^*;H^1(\mathbb{R}_{\pm}))\Big\}.
\ea
$$

In order to prove Proposition \ref{pro-1} for our purpose, we shall combine a local existence result in Proposition \ref{local},  together with a-priori estimates in Proposition \ref{priori} by continuation arguments. To state the local existence result precisely, we set
$$
(\bar\phi,\bar\psi,\bar\z)(\tau,y)=(v,u,\t)(\tau,y)-(V,U,\T)(y-s_3\tau),
$$
then $(\bar\phi,\bar\psi,\bar\z)$ satisfies \eqref{perturb} without time shift function $\mb{X}(\tau)$, and has the same initial value with $(\phi_0,\psi_0,\z_0)$ since $\mb{X}(0)=0$. The local existence of a small solution $(\bar\phi,\bar\psi,\bar\z)$ has been established by \cite{Hoff}, we just state it and omit its proof for brevity.

\begin{proposition}\label{local}
{\rm (Local existence)} For any $\Xi>0$, suppose the initial data and the wave strength satisfy $\mathcal{N}_0+\d\leq \Xi$. Then there exists a positive constant time $\tau_0$ depending on $\Xi$ such that problem \eqref{perturb}, \eqref{per-initial} with $\mb{X}(\tau)=0$ has a unique solution $(\bar\phi,\bar\psi,\bar\z)(\tau,y)\in \mathcal{X}(0,\tau_0)$ satisfying
$$
A(\tau_0)+B(\tau_0)+F(\tau_0)\leq 4(\mathcal{N}_0+\d)^2,
$$
where
$$
\ba
A(\tau_0)&=\sup_{0\leq \tau\leq \tau_0}\Big(\|(\bar\phi,\bar\psi,\bar\z)(\tau,\cdot)\|^2+\|\bar\phi_y(\tau,\cdot)\rlap{--}\hspace{0.07em}\|^2
\Big)
+\int_0^{\tau_0}\|(\bar\psi_y,\bar\z_y)\|^2d\tau,\\
B(\tau_0)&=\sup_{0\leq \tau\leq \tau_0}\Big(g(\tau)^{\f12}\|\bar\psi_y(\tau,\cdot)\|^2+g(\tau)\|\bar\z_y(\tau,\cdot)\|^2\Big)
+\int_{0}^{\tau_0}g(\tau)^{\f12+\vartheta}\Big(\|\bar\psi_\tau\|^2+\Big\|\Big(\f{\bar\psi_y}{v}\Big)_y\rlap{--}\hspace{0.03em}\Big\|^2\Big)d\tau\\
&+\int_0^{\tau_0}g(\tau)\Big(\|\bar\psi_y^2\|^2+\|\bar\z_\tau\|^2+\Big\|\Big(\f{\bar\z_y}{v}\Big)_y\rlap{--}\hspace{0.03em}\Big\|^2\Big)d\tau,\\
F(\tau_0)&=\sup_{0\leq \tau\leq \tau_0}
\Big\{g(\tau)^{\f32+\vartheta}\Big(\|\bar\psi_\tau(\tau,\cdot)\|^2+\Big\|\Big(\f{\bar\psi_y}{v}\Big)_y(\tau,\cdot)\rlap{--}\hspace{0.03em}\Big\|^2\Big)
+g(\tau)^{3}\Big(\|\bar\z_\tau(\tau,\cdot)\|^2+\Big\|\Big(\f{\bar\z_y}{v}\Big)_y(\tau,\cdot)\rlap{--}\hspace{0.03em}\Big\|^2\Big)
\Big\}\\
&+\int_0^{\tau_0}\Big(g(\tau)^{\f32+\vartheta}\|\bar\psi_{\tau y}\rlap{--}\|^2
+g(\tau)^3\|\bar\z_{\tau y}\rlap{--}\|^2\Big)d\tau,
\ea
$$
where $g(\tau)=\tau\wedge1=\min\{\tau,1\}$ and $\vartheta\in(0,1)$. Moreover, $v$, $u$, $\theta$ have the same regularity as in Theorem \ref{thm}. Thus, $v$, $u_y$, $\t_y$ have one-side limit at $y=0$ and satisfy the jump conditions
$$
\Big[p-\f{u_y}{v}\Big]=\Big[\f{\t_y}{v}\Big]=0.
$$
Finally, one has the following estimate on the jump at $y=0$,
$$
|[v](\tau)|\leq C\d e^{-c\tau},\quad \tau>0
$$
for some positive constants $C$ and $c$ independent of time $\tau$.
\end{proposition}

Next, we show the following (uniform) a-priori estimates.

\begin{proposition}\label{priori}
{\rm (A-priori estimates)} Let $(\phi,\psi,\z)\in\mathcal{X}[\tau_1,\tau_2]$ be a solution to the problem \eqref{perturb}, \eqref{per-initial} for $0<\tau_1\leq \tau_2$. Then there exist  positive constants $\chi_0\leqq1$ and $C_0$ independent of $\tau_1$, $\tau_2$, such that if
\begin{equation}\label{assumption}
\ba
N(\tau_1,\tau_2)+\d\leq \chi_0,
\ea
\end{equation}
then it holds
\begin{equation}\label{full-es}
\ba
&N^2(\tau_1,\tau_2)
+\int_{\tau_1}^{\tau_2}\left(\d|\dot{\mathbf{X}}(\tau)|^2+\Big\|\sqrt{|U_y^{-\mathbf{X}}|}\,\psi\Big\|^2
+\|\phi_y\rlap{--}\hspace{0.07em}\|^2+\|(\psi_y,\z_y)\rlap{--}\hspace{0.07em}\|_{H^1}^2
+\|(\psi_{\tau y},\z_{\tau y})\rlap{--}\hspace{0.07em}\|^2\right)d\tau\\
&\leq C_0(N^2(\tau_1)+\d).
\ea
\end{equation}
In addition, if the shift function $\mathbf{X}(\tau)$ is an absolutely continuous solution to \eqref{X}, then it holds for $\tau\in[\tau_1,\tau_2]$,
\be\label{point-X}
|\dot{ \mathbf{X}}(\tau)|\leq C_0\|\psi(\tau,\cdot)\|_{L^{\infty}}.
\ee
\end{proposition}

The proof of Proposition \ref{priori} will be given in Sect. \ref{pf-priori}.
Then we can demonstrate Proposition \ref{pro-1} by the continuation arguments based on Propositions \ref{local} and \ref{priori}, the detailed proof can be found in \cite{WW}.

%
%
%
%

\section{A-priori estimates}\label{pf-priori}
\setcounter{equation}{0}

Throughout this section, $C$ denotes a positive constant which may change from
line to line, but which stays independent of the wave strength $\d$. From now on until the end of this paper, we always assume that $\chi_0\leqq1$.

\subsection{Relative entropy method}

The following lemma is derived by the relative entropy method.

\begin{lemma}\label{le-1}
Let $a(\tau,y)$ be the weighted function defined by \eqref{a}, and $\mb{X}(\tau)$ be the shift defined by \eqref{X}. Then it holds
\be\label{ba-expression}
\f{d}{d\tau}\int_{\mathbb{R}}\hspace{-1.15em}-
a^{-\mathbf{X}}\mathcal{E}dy+\mathbf{G}_1(\tau)+\mathbf{G}_2(\tau)
+\mathbf{D}_1(\tau)+\mathbf{D}_2(\tau)+\f{\d}{M}|\dot{\mathbf{X}}(\tau)|^2\\
=\dot{\mathbf{X}}(\tau)\sum_{i=1}^3\mathbf{Y}_{i}(\tau)+\sum_{i=1}^{7}\mathbf{B}_{i}(\tau),
\ee
where $\mb{X}(\tau)$ has been defined in \eqref{X}, and
$$
\mathcal{E}=R\T^{-\mb{X}}\Phi\Big(\f{v}{V^{-\mb{X}}}\Big)+\f{\psi^2}{2}+\f{R \T^{-\mb{X}}}{\gamma-1}\Phi\Big(\f{\t}{\T^{-\mb{X}}}\Big),\quad
\Phi(z)=z-\ln z-1,
$$
$$
\ba
\mathbf{G}_{1}(\tau)&=\int_{\mathbb{R}}\hspace{-1.15em}-
a^{-\mathbf{X}}_y\f{s_3p_+}{2v_+}\left(\phi+\f{\psi}{s_3}\right)^2dy,\\
\mathbf{G}_{2}(\tau)&=\int_{\mathbb{R}}\hspace{-1.15em}-
a^{-\mathbf{X}}_y\f{R}{\gamma-1}\f{s_3}{2\t_+}\left(\z-\f{\gamma-1}{R}\f{p_+\psi}{s_3}\right)^2dy,\\
\mathbf{D}_1(\tau)&=\int_{\mathbb{R}}
a^{-\mathbf{X}}\f{\T^{-\mb{X}}}{v\t}\psi_y^2dy, \qquad \mathbf{D}_2(\tau)=\nu\int_{\mathbb{R}}
a^{-\mathbf{X}}\f{\T^{-\mb{X}}}{v\t^2}\z_y^2dy,
\ea
$$
\begin{align} \nm
\mathbf{Y}_{1}(\tau)&=\int_{\mathbb{R}}\hspace{-1.15em}-
a^{-\mathbf{X}}\left(
\f{P^{-\mb{X}}}{V^{-\mb{X}}}\Big(\phi+\f{\psi}{s_3}\Big)V_y^{-\mb{X}}
+\f{R}{\gamma-1}\Big(\z-\f{\gamma-1}{R}\f{p_+\psi}{s_3}\Big)\f{\T_y^{-\mb{X}}}{\T^{-\mb{X}}}\right) dy,\\  \nm
\mathbf{Y}_{2}(\tau)&=-\int_{\mathbb{R}}\hspace{-1.15em}-
a^{-\mathbf{X}}\T_y^{-\mathbf{X}}\left(R\Phi\Big(\f{v}{V^{-\mb{X}}}\Big)
+\f{R}{\gamma-1}\Phi\Big(\f{\t}{\T^{-\mb{X}}}\Big)
\right)dy,\\ \nm
\mathbf{Y}_{3}(\tau)&=-\int_{\mathbb{R}}\hspace{-1.15em}-a^{-\mathbf{X}}_y\mathcal{E}dy,
\end{align}
and
\begin{align} \nm
\mathbf{B}_{1}(\tau)&=\int_{\mathbb{R}}\hspace{-1.15em}-
-a^{-\mathbf{X}}U^{-\mathbf{X}}_y P^{-\mb{X}}\left(\Phi\Big(\f{V^{-\mb{X}}\t}{\T^{-\mb{X}} v}\Big)+\gamma\Phi\Big(\f{v}{V^{-\mb{X}}}\Big)\right)dy
\\ \nm
&+\int_{\mathbb{R}}\hspace{-1.15em}-
a^{-\mathbf{X}}\left(\nu\Big(\f{\T_y^{-\mathbf{X}}}{V^{-\mathbf{X}}}\Big)_y
+\f{(U_y^{-\mathbf{X}})^2}{V^{-\mathbf{X}}}\right)\left((\gamma-1)\Phi\Big(\f{v}{ V^{-\mb{X}}}\Big)-\Phi\Big(\f{\T^{-\mb{X}}}{\t }\Big)\right) dy\\ \nm
&+\int_{\mathbb{R}}\hspace{-1.15em}-
-\f{a^{-\mathbf{X}}_y}{2s_3}\left(s_3^2-\f{\gamma p_+}{v_+}\right)\psi^2dy=:\sum_{i=1}^3\mb{B}_{1,i}(t),\\ \nm
\mathbf{B}_{2}(\tau)&=O(1)(N(\tau_1,\tau_2)+\d)\int_{\mathbb{R}}\hspace{-1.25em}-
a^{-\mathbf{X}}_y
\left|\left(\phi+\f{\psi}{s_3}, \z-\f{\gamma-1}{R}\f{p_+\psi}{s_3}\right)\right|^2dy,\\ \nm
\mathbf{B}_{3}(\tau)&=O(1)\int_{\mathbb{R}}\hspace{-1.15em}-
a^{-\mathbf{X}}_y
\big(|(\phi,\z)|+\d\big)\psi^2dy,\\ \nm
\mathbf{B}_{4}(\tau)&=\int_{\mathbb{R}}\hspace{-1.15em}-
a^{-\mathbf{X}}
\left(\f{ U^{-\mb{X}}_y\phi}{V^{-\mb{X}} v}\psi_y
+2\f{U^{-\mb{X}}_y\z}{v\t}\psi_y+\nu\f{\T^{-\mb{X}}_y\z}{v\t^2}\z_y
+\nu\f{\T^{-\mb{X}}}{V^{-\mb{X}}}\f{\T^{-\mb{X}}_y\phi}{v\t^2}\z_y\right)dy,\\ \nm
\mathbf{B}_{5}(\tau)&=-\int_{\mathbb{R}}\hspace{-1.15em}-
a^{-\mathbf{X}}
\left(\f{(U^{-\mb{X}}_y)^2}{V^{-\mb{X}}}\f{\phi\z}{v\t}
+\nu\f{(\T^{-\mb{X}}_y)^2}{V^{-\mb{X}}}\f{\phi\z}{v\t^2}\right)dy,\\ \nm
\mathbf{B}_{6}(\tau)&=-\int_{\mathbb{R}}\hspace{-1.15em}-
 a^{-\mathbf{X}}_y\left(\f{\psi\psi_y}{v}+\nu\f{\z\z_y}{v\t}\right)
 dy,\\ \nm
\mathbf{B}_{7}(\tau)&=\int_{\mathbb{R}}\hspace{-1.15em}-
a^{-\mathbf{X}}_y
\left(\f{U^{-\mb{X}}_y}{V^{-\mb{X}}}\f{\phi\psi}{v}+\nu\f{\T^{-\mb{X}}_y}{V^{-\mb{X}}}\f{\phi\z}{v\t}\right)dy.
\end{align}

\end{lemma}

\begin{remark}
Since $s_3 a_y^{-\mathbf{X}}>0$, thus $\mathbf{G}_{i}(\tau)(i=1,2)$ consists of terms with good sign, while $\mathbf{B}_{i}(\tau)(i=1,\cdots,7)$ consists of bad terms.
\end{remark}

{\textbf{\emph{Proof}}:} It follows from \eqref{a} that
\be\label{a-x}
\ba
a^{-\mathbf{X}}:&=1+\f{p_- -P^{-\mathbf{X}}}{\sqrt{\d}}
=1+\f{p_- -P(V^{-\mathbf{X}},\T^{-\mathbf{X}})}{\sqrt{\d}}\\
&=1+\f{p_- -P(V,\T)(y-s_3\tau-\mathbf{X}(\tau))}{\sqrt{\d}}.
\ea
\ee
Then it holds
$$
|a^{-\mathbf{X}}-1|\leq \sqrt{\d},\qquad a_\tau^{-\mathbf{X}}=-s_3a_y^{-\mathbf{X}}-\dot{\mathbf{X}}(\tau)a_y^{-\mathbf{X}}.
$$
First, we multiply \eqref{perturb}$_1$ by $-a^{-\mathbf{X}}R\T^{-\mb{X}}(1/v-1/V^{-\mb{X}})$ to have
\be\label{ba-1}
\ba
-a^{-\mathbf{X}}R\T^{-\mb{X}}\left(\f{1}{v}-\f{1}{ V^{-\mb{X}}}\right)\phi_\tau
&+\dot{\mathbf{X}}(\tau)a^{-\mathbf{X}}R\T^{-\mb{X}}\left(\f{1}{v}-\f{1}{ V^{-\mb{X}}}\right)V_y^{-\mathbf{X}}\\
&+a^{-\mathbf{X}}R\T^{-\mb{X}}\left(\f{1}{v}-\f{1}{ V^{-\mb{X}}}\right)\psi_y=0,
\ea
\ee
where
$$
\ba
-a^{-\mathbf{X}}R\T^{-\mb{X}}\left(\f{1}{v}-\f{1}{V^{-\mb{X}}}\right)\phi_\tau
=\left(a^{-\mathbf{X}}R\T^{-\mb{X}}\Phi\Big(\f{v}{V^{-\mb{X}}}\Big)\right)_\tau
+s_3a_y^{-\mathbf{X}}R\T^{-\mb{X}}\Phi\Big(\f{v}{V^{-\mb{X}}}\Big)\\
+\dot{\mathbf{X}}(\tau)a_y^{-\mathbf{X}}R\T^{-\mb{X}}\Phi\Big(\f{v}{V^{-\mb{X}}}\Big)
+a^{-\mathbf{X}}V_\tau^{-\mb{X}}\f{P^{-\mb{X}}\phi^2}{V^{-\mb{X}} v}-a^{-\mathbf{X}}R\T^{-\mb{X}}_\tau\Phi\Big(\f{v}{V^{-\mb{X}}}\Big).
\ea
$$
Next, we multiply \eqref{perturb}$_2$ by $a^{-\mathbf{X}}\psi$ to have
\begin{align}\label{ba-2} \nm
&\quad\left(a^{-\mathbf{X}}\f{\psi^2}{2}\right)_\tau+s_3a^{-\mathbf{X}}_y\f{\psi^2}{2}
+\dot{\mathbf{X}}(\tau)a^{-\mathbf{X}}_y\f{\psi^2}{2}-\dot{\mathbf{X}}(\tau)a^{-\mathbf{X}}\psi U_y^{-\mathbf{X}}\\ \nm
&\quad+\left(a^{-\mathbf{X}}(p-P^{-\mb{X}})\psi\right)_y-a^{-\mathbf{X}}(p-P^{-\mb{X}} )\psi_y-a^{-\mathbf{X}}_y(p-P^{-\mb{X}})\psi\\
&=\left(a^{-\mathbf{X}}\Big(\f{u_y}{v}-\f{U^{-\mb{X}}_y}{V^{-\mb{X}}}\Big)\psi\right)_y
- a^{-\mathbf{X}}\f{\psi_y^2}{v}
- a^{-\mathbf{X}}U^{-\mb{X}}_y\left(\f{1}{v}-\f{1}{V^{-\mb{X}}}\right)\psi_y\\ \nm
&- a^{-\mathbf{X}}_y\f{\psi\psi_y}{v}- a_y^{-\mathbf{X}} U^{-\mb{X}}_y\left(\f{1}{v}-\f{1}{V^{-\mb{X}}}\right)\psi.
\end{align}
Third, we multiply \eqref{perturb}$_3$ by $a^{-\mathbf{X}}\z/\t$ to have
\begin{align} \label{ba-3} \nm
&\quad a^{-\mathbf{X}}\f{R}{\gamma-1}\f{\z \z_\tau}{\t}-\dot{\mathbf{X}}(\tau)a^{-\mathbf{X}}\f{\z}{\t}\f{R\T_y^{-\mb{X}}}{\gamma-1}
+a^{-\mathbf{X}}\f{R\z}{v}\psi_y+a^{-\mathbf{X}}(p-P^{-\mb{X}}) U^{-\mb{X}}_y\f{\z}{\t}\\
&=\nu\left(a^{-\mathbf{X}}\Big(\f{\t_y}{v}-\f{\T^{-\mb{X}}_y}{V^{-\mb{X}} }\Big)\f{\z}{\t}\right)_y
-\nu a^{-\mathbf{X}}\f{\T^{-\mb{X}}\z_y^2}{v\t^2}+\nu a^{-\mathbf{X}}\f{\T^{-\mb{X}}_y\z}{v\t^2}\z_y
\\ \nm
&-\nu a^{-\mathbf{X}}\f{\T^{-\mb{X}}\T^{-\mb{X}}_y}{\t^2}\left(\f{1}{v}-\f{1}{V^{-\mb{X}} }\right)\z_y
+\nu a^{-\mathbf{X}}\f{(\T^{-\mb{X}}_y)^2}{\t^2}\left(\f{1}{v}-\f{1}{V^{-\mb{X}}}\right)\z
-\nu a^{-\mathbf{X}}_y\f{\z\z_y}{v\t}\\ \nm
&-\nu a_y^{-\mathbf{X}}\T^{-\mb{X}}_y \left(\f{1}{v}-\f{1}{V^{-\mb{X}}}\right)\f{\z}{\t}
+ a^{-\mathbf{X}}\f{\psi_y^2\z}{v\t}
+2 a^{-\mathbf{X}}\f{\psi_y U^{-\mb{X}}_y}{v}\f{\z}{\t}+ a^{-\mathbf{X}} (U^{-\mb{X}}_y)^2\left(\f{1}{v}-\f{1}{V^{-\mb{X}}}\right)\f{\z}{\t},
\end{align}
where
$$
\ba
a^{-\mathbf{X}}\f{R}{\gamma-1}\f{\z \z_\tau}{\t}
&=\left(a^{-\mathbf{X}}\f{R\T^{-\mb{X}}}{\gamma-1}\Phi\Big(\f{\t}{\T^{-\mb{X}}}\Big)\right)_\tau
+s_3a_y^{-\mathbf{X}}\f{R\T^{-\mb{X}}}{\gamma-1}\Phi\Big(\f{\t}{\T^{-\mb{X}}}\Big)\\
&+\dot{\mathbf{X}}(\tau)a_y^{-\mathbf{X}}\f{R\T^{-\mb{X}}}{\gamma-1}\Phi\Big(\f{\t}{\T^{-\mb{X}}}\Big)
+a^{-\mathbf{X}}\f{R\T^{-\mb{X}}_\tau}{\gamma-1}\Phi\Big(\f{\T^{-\mb{X}}}{\t}\Big).
\ea
$$

Direct calculation yields
$$
\ba
&\quad V^{-\mb{X}}_\tau\f{P^{-\mb{X}}\phi^2}{V^{-\mb{X}} v}-R\T^{-\mb{X}}_\tau\Phi\Big(\f{v}{V^{-\mb{X}}}\Big)
+\f{R\T^{-\mb{X}}_\tau}{\gamma-1}\Phi\Big(\f{\T^{-\mb{X}}}{\t}\Big)+(p-P^{-\mb{X}}) U^{-\mb{X}}_y\f{\z}{\t}\\
&=U^{-\mathbf{X}}_y P^{-\mb{X}}\left(\Phi\Big(\f{V^{-\mb{X}}\t}{\T^{-\mb{X}} v}\Big)+\gamma\Phi\Big(\f{v}{V^{-\mb{X}}}\Big)\right)\\ \nm
&-\left(\nu\Big(\f{\T_y^{-\mathbf{X}}}{V^{-\mathbf{X}}}\Big)_y
+\f{(U_y^{-\mathbf{X}})^2}{V^{-\mathbf{X}}}\right)\left((\gamma-1)\Phi\Big(\f{v}{ V^{-\mb{X}}}\Big)
-\Phi\Big(\f{\T^{-\mb{X}}}{\t }\Big)\right)\\ \nm
&-\dot{\mathbf{X}}(\tau)\left(\f{P^{-\mb{X}}\phi^2}{V^{-\mb{X}}}V_y^{-\mathbf{X}}
-\Phi\Big(\f{v}{V^{-\mb{X}}}\Big)R\T_y^{-\mathbf{X}}+\Phi\Big(\f{\T^{-\mb{X}}}{\t}\Big)
\f{R\T_y^{-\mathbf{X}}}{\gamma-1}\right).
\ea
$$

We add \eqref{ba-1}, \eqref{ba-2}, \eqref{ba-3} and use the above equality to have
\begin{align} \nm
&(a^{-\mathbf{X}}\mathcal{E})_\tau+a^{-\mathbf{X}}_y\big(s_3\mathcal{E}-(p-P^{-\mb{X}} )\psi\big)
+\big(a^{-\mathbf{X}}(p-P^{-\mb{X}})\psi\big)_y\\ \nm
&+\dot{\mathbf{X}}(\tau)a^{-\mathbf{X}}_y\mathcal{E}
-\dot{\mathbf{X}}(\tau)a^{-\mathbf{X}}\left(
\f{P^{-\mb{X}}}{v}\phi V^{-\mathbf{X}}_y+\psi U^{-\mathbf{X}}_y
+\f{\z}{\t}\f{R\T_y^{-\mathbf{X}}}{\gamma-1}
\right)\\ \nm
&-\dot{\mathbf{X}}(\tau)a^{-\mathbf{X}}\left(
\f{P^{-\mb{X}}\phi^2}{V^{-\mb{X}}v}V^{-\mathbf{X}}_y
-\Phi\Big(\f{v}{V^{-\mb{X}}}\Big)R\T_y^{-\mathbf{X}}
+\Phi\Big(\f{\T^{-\mb{X}}}{\t}\Big)\f{R\T_y^{-\mathbf{X}}}{\gamma-1}
\right)\\ \nm
&+a^{-\mathbf{X}}U^{-\mathbf{X}}_y P^{-\mb{X}}\left(\Phi\Big(\f{V^{-\mb{X}}\t}{\T^{-\mb{X}} v}\Big)+\gamma\Phi\Big(\f{v}{V^{-\mb{X}}}\Big)\right)\\ \nm
&-a^{-\mathbf{X}}\left(\nu\Big(\f{\T_y^{-\mathbf{X}}}{V^{-\mathbf{X}}}\Big)_y
+\f{(U_y^{-\mathbf{X}})^2}{V^{-\mathbf{X}}}\right)\left((\gamma-1)\Phi\Big(\f{v}{V^{-\mb{X}} }\Big)
-\Phi\Big(\f{\T^{-\mb{X}}}{\t }\Big)\right)\\ \label{L2-diff}
&=\left(a^{-\mathbf{X}}
\Big(\f{u_y}{v}-\f{U^{-\mb{X}}_y}{V^{-\mb{X}}}\Big)\psi
+\nu a^{-\mathbf{X}}\Big(\f{\t_y}{v}-\f{\T^{-\mb{X}}_y}{V^{-\mb{X}}}\Big)\f{\z}{\t}
\right)_y
- a^{-\mathbf{X}}\f{\T^{-\mb{X}}\psi_y^2}{v\t}-\nu a^{-\mathbf{X}}\f{\T^{-\mb{X}}\z_y^2}{v\t^2}\\ \nm
& +a^{-\mathbf{X}}\left(\f{U^{-\mb{X}}_y\phi}{V^{-\mb{X}}v}\psi_y
+2\f{U^{-\mb{X}}_y\z}{v\t}\psi_y+\nu\f{\T^{-\mb{X}}_y\z}{v\t^2}\z_y
+\nu\f{\T^{-\mb{X}}}{V^{-\mb{X}}}\f{\T^{-\mb{X}}_y\phi}{v\t^2}\z_y\right)\\ \nm
& -a^{-\mathbf{X}}\left(\f{ (U^{-\mb{X}}_y)^2}{V^{-\mb{X}}}\f{\phi\z}{v\t}
+\nu\f{(\T^{-\mb{X}}_y)^2}{V^{-\mb{X}}}\f{\phi\z}{v\t^2}\right)
-a^{-\mathbf{X}}_y\left(\f{\psi\psi_y}{v}+\nu\f{\z\z_y}{v\t}\right)\\ \nm
&+a^{-\mathbf{X}}_y\left(\f{ U^{-\mb{X}}_y}{V^{-\mb{X}}}\f{\phi\psi}{v}+\nu\f{\T^{-\mb{X}}_y}{V^{-\mb{X}} }\f{\phi\z}{v\t}\right).
\end{align}
We directly calculate to have
\begin{align}\nm
&\quad\, s_3\mathcal{E}-(p-P^{-\mb{X}})\psi \\ \nm
&=s_3\left(R\T^{-\mb{X}}\Phi\Big(\f{v}{V^{-\mb{X}}}\Big)+\f{\psi^2}{2}
+\f{R\T^{-\mb{X}}}{\gamma-1}\Phi\Big(\f{\t}{\T^{-\mb{X}}}\Big) \right)-R\left(\f{\t}{v}-\f{\T^{-\mb{X}}}{V^{-\mb{X}}}\right)\psi\\ \label{s3-eta}
&=\f{s_3p_+}{2v_+}\left(\phi+\f{\psi}{s_3}\right)^2
+\f{R}{\gamma-1}\f{s_3}{2\t_+}\left(\z-\f{\gamma-1}{R}\f{p_+\psi}{s_3}\right)^2
+\left(s_3^2-\f{\gamma p_+}{v_+}\right)\f{\psi^2}{2s_3}\\ \nm
&\quad+O(1)(N(\tau_1,\tau_2)+\d)\left|\left(\phi+\f{\psi}{s_3}, \z-\f{\gamma-1}{R}\f{p_+\psi}{s_3}\right)\right|^2+O(1)\big(|(\phi,\z)|+\d\big)\psi^2,
\end{align}
and
\begin{align}\nm
&-\dot{\mathbf{X}}(\tau)a^{-\mathbf{X}}\left(
\f{P^{-\mb{X}}}{v}\phi V^{-\mathbf{X}}_y+\psi U^{-\mathbf{X}}_y
+\f{\z}{\t}\f{R \T_y^{-\mathbf{X}}}{\gamma-1}
\right)\\ \label{simplify}
&-\dot{\mathbf{X}}(\tau)a^{-\mathbf{X}}\left(
\f{P^{-\mb{X}}\phi^2}{V^{-\mb{X}} v}V^{-\mathbf{X}}_y
-\Phi\Big(\f{v}{V^{-\mb{X}}}\Big)R \T_y^{-\mathbf{X}}
+\Phi\Big(\f{\T^{-\mb{X}}}{\t}\Big)\f{R \T_y^{-\mathbf{X}}}{\gamma-1}
\right)\\ \nm
=&-\dot{\mathbf{X}}(\tau)a^{-\mathbf{X}}\left(
\f{P^{-\mb{X}}}{V^{-\mb{X}}}\phi V^{-\mathbf{X}}_y+\psi U^{-\mathbf{X}}_y
+\f{\z}{\T^{-\mb{X}}}\f{R \T_y^{-\mathbf{X}}}{\gamma-1}
\right)\\ \nm
&+\dot{\mathbf{X}}(\tau)a^{-\mathbf{X}}\T_y^{-\mathbf{X}}
\left(R\Phi\Big(\f{v}{V^{-\mb{X}}}\Big)+\f{R}{\gamma-1}\Phi\Big(\f{\t}{\T^{-\mb{X}}}\Big)
\right).
\end{align}
Integrating \eqref{L2-diff} with respect to $y$ over $\mathbb{R}_{\pm}$, the jump term in \eqref{L2-diff} across $y=0$ vanishes, i.e.,
$$
\ba
&\left[a^{-\mathbf{X}}(p-P^{-\mb{X}})\psi-a^{-\mathbf{X}}
\Big(\f{u_y}{v}-\f{U^{-\mb{X}}_y}{V^{-\mb{X}}}\Big)\psi
-\nu a^{-\mathbf{X}}\Big(\f{\t_y}{v}-\f{\T^{-\mb{X}}_y}{V^{-\mb{X}}}\Big)\f{\z}{\t}\right](\tau)\\
=&\left[a^{-\mb{X}}\Big(p-\f{u_y}{v}\Big)\psi\right](\tau)
-\left[a^{-\mb{X}}\Big(P-\f{U^{-\mb{X}}_y}{V^{-\mb{X}}}\Big)\psi\right](\tau)
-\nu \left[a^{-\mb{X}}\Big(\f{\t_y}{v}-\f{\T^{-\mb{X}}_y}{V^{-\mb{X}}}\Big)\f{\z}{\t}\right](\tau)\\
=&\left[p-\f{u_y}{v}\right]a^{-\mb{X}}\psi(\tau,0)
-\left[P-\f{U^{-\mb{X}}_y}{V^{-\mb{X}}}\right]a^{-\mb{X}}\psi(\tau,0)
-\nu \left(\Big[\f{\t_y}{v}\Big]-\Big[\f{\T^{-\mb{X}}_y}{V^{-\mb{X}}}\Big]\right)a^{-\mb{X}}\f{\z(\tau,0)}{\t(\tau,0)}=0.
\ea
$$
Using the above equality and \eqref{s3-eta}, \eqref{simplify}, we have
\be\label{ba-expression-1}
\f{d}{d\tau}\int_{\mathbb{R}}\hspace{-1.15em}-
a^{-\mathbf{X}}\mathcal{E}dy+\mathbf{G}_1(\tau)+\mathbf{G}_2(\tau)
+\mathbf{D}_1(\tau)+\mathbf{D}_2(\tau)
=\dot{\mathbf{X}}(\tau)\mathbf{Y}(\tau)+\sum_{i=1}^{7}\mathbf{B}_{i}(\tau),
\ee
where
$$
\ba
\mathbf{Y}(\tau)&=\int_{\mathbb{R}}\hspace{-1.15em}-a^{-\mathbf{X}}\left(
\f{P^{-\mb{X}}}{V^{-\mb{X}}}\phi V^{-\mathbf{X}}_y+\psi U^{-\mathbf{X}}_y
+\f{\z}{\T^{-\mb{X}}}\f{R \T_y^{-\mathbf{X}}}{\gamma-1}
\right)dy\\
&-\int_{\mathbb{R}}\hspace{-1.15em}-
a^{-\mathbf{X}}\T_y^{-\mathbf{X}}\left(R\Phi\Big(\f{v}{V^{-\mb{X}} }\Big)
+\f{R}{\gamma-1}\Phi\Big(\f{\t}{\T^{-\mb{X}}}\Big)\right)dy
-\int_{\mathbb{R}}\hspace{-1.15em}-a^{-\mathbf{X}}_y\mathcal{E}dy\\
&=-\f{\d}{M}\dot{\mb{X}}(\tau)+\sum_{i=1}^3\mathbf{Y}_{i}(\tau),
\ea
$$
by using the definition of $\dot{\mb{X}}(\tau)$ in \eqref{X} and the definition of $\mb{Y}_i(\tau)(i=1,2,3)$. Thus, it holds
\be\label{XY}
\dot{\mathbf{X}}(\tau)\mathbf{Y}(\tau)=-\f{\d}{M}|\dot{\mathbf{X}}(\tau)|^2
+\dot{\mathbf{X}}(\tau)\sum_{i=1}^3\mathbf{Y}_{i}(\tau).
\ee
Combination of \eqref{ba-expression-1} and \eqref{XY} gives \eqref{ba-expression}. The proof of Lemma \ref{le-1} is completed.

\hfill $\Box$

%
%
%
%

\subsection{Contraction property}
The contraction property is expressed in the following Lemma.

\begin{lemma}\label{le-contraction}
Under the hypotheses of Proposition \ref{priori}, there exists constant $C>0$, such that for any $\tau\in[\tau_1,\tau_2]$, it holds
\be\label{es-contaction}
\mathbf{B}_1(\tau)-\f78\big(
\mathbf{D}_1(\tau)+\mathbf{D}_2(\tau)\big)-\f{\d}{2M}|\dot{\mathbf{X}}(\tau)|^2
\leq -\f{\alpha_+}{3}\int_{\mathbb{R}}|P^{-\mathbf{X}}_y|\psi^2dy
+C\d^{\f14}\big(\mathbf{G}_1(\tau)+\mathbf{G}_2(\tau)\big),
\ee
where $\alpha_+:=\f{\gamma+1}{2\gamma}\f{s_3}{p_+}$.
\end{lemma}

{\textbf{\emph{Proof}}:} We introduce a new variable
\be\label{z}
z=\f{p_- -P^{-\mathbf{X}}}{\d},
\ee
which maps $y\in\mathbb{R}$ to $z\in[0,1]$. Then it follows from \eqref{a-x} and \eqref{z} that
\be\label{change-variable}
\f{dz}{dy}=-\f{1}{\d}P^{-\mathbf{X}}_y=\f{1}{\sqrt{\d}}a^{-\mathbf{X}}_y.
\ee

We can deduce from the viscous shock equations \eqref{shock-equ} that
\be\label{equivalent}
\left\{
\ba
V^{-\mathbf{X}}_y&=-\f{1}{s_3^2}P^{-\mathbf{X}}_y+R^s_1,\\
U^{-\mathbf{X}}_y&=\f{1}{s_3}P^{-\mathbf{X}}_y+R^s_2,\\
\T^{-\mathbf{X}}_y&=\f{\gamma-1}{R s_3^2}P^{-\mathbf{X}}P^{-\mathbf{X}}_y+R^s_3,
\ea
\right.
\ee
where $|R^s_i|\leq C\d|P^{-\mathbf{X}}_y|~(i=1,2,3)$.


$\bullet$ \textbf{Change of variable for $\mathbf{B}_1(\tau)$}: We first treat each term $\mb{B}_{1,i}(\tau)~(i=1,\cdots,3)$ by using a new variable \eqref{z}. By Taylor expansion and  \eqref{equivalent}, it holds
$$
\ba
\mb{B}_{1,1}(\tau)&\leq\int_{\mathbb{R}}\hspace{-1.15em}-
-U^{-\mathbf{X}}_y P^{-\mb{X}}\Big[\f12\Big(\f{V^{-\mb{X}}\t}{\T^{-\mb{X}}v}-1\Big)^2
+\f{\gamma}{2}\Big(\f{v}{V^{-\mb{X}}}-1\Big)^2\Big]dy
+C(N(\tau_1,\tau_2)+\sqrt{\d})\int_{\mathbb{R}}\hspace{-1.15em}-|U^{-\mathbf{X}}_y|\,|(\phi,\z)|^2dy\\
&\leq\int_{\mathbb{R}}\hspace{-1.15em}-
-U^{-\mathbf{X}}_y p_+\Big[\f12\Big(\f{\z}{\t_+}-\f{\phi}{v_+}\Big)^2+\f{\gamma}{2}\f{\phi^2}{v_+^2}\Big]dy
+C(N(\tau_1,\tau_2)+\sqrt{\d})\int_{\mathbb{R}}\hspace{-1.15em}-|U^{-\mathbf{X}}_y|\,|(\phi,\z)|^2dy\\
&\leq\int_{\mathbb{R}}\hspace{-1.15em}-
-P_y^{-\mb{X}}\f{p_+}{s_3}\Big(\f{\gamma+1}{2}\f{\phi^2}{v_+^2}+\f{1}{2}\f{\z^2}{\t_+^2}
-\f{\phi\z}{v_+\t_+}\Big)dy
+C(N(\tau_1,\tau_2)+\sqrt{\d})\int_{\mathbb{R}}\hspace{-1.15em}-|P^{-\mathbf{X}}_y|\,|(\phi,\z)|^2dy.
\ea
$$
Using inequality $m^2\leq (1+\eta)n^2+(1+1/\eta)(m-n)^2$ for $\eta>0$,  and \eqref{equivalent}, it holds
$$
\ba
&\quad\int_{\mathbb{R}}\hspace{-1.15em}-
-P_y^{-\mb{X}}\f{p_+}{s_3}\Big(\f{\gamma+1}{2}\f{\phi^2}{v_+^2}+\f{1}{2}\f{\z^2}{\t_+^2}\Big)dy\\
&\leq (1+\eta)\int_{\mathbb{R}}
-P^{-\mathbf{X}}_y \f{p_+}{s_3}
\left(\f{\gamma+1}{2v_+^2}\f{\psi^2}{s_3^2}
+\f{(\gamma-1)^2}{2R^2}\f{p_+^2\psi^2}{\t_+^2s_3^2}
\right)dy\\
&\quad+\f{C}{\eta}\int_{\mathbb{R}}\hspace{-1.15em}-|P^{-\mathbf{X}}_y|\left\{
\Big(\phi+\f{\psi}{s_3}\Big)^2+\Big(\z-\f{\gamma-1}{R}\f{p_+\psi}{s_3}\Big)^2\right\}dy\\
&=\int_{\mathbb{R}}\hspace{-1.15em}-
-P_y^{-\mb{X}}\f{\gamma^2-\gamma+2}{2}\f{p_+\psi^2}{s_3^3v_+^2}dy
+\f{C}{\eta}\sqrt{\d}(\mb{G}_1(\tau)+\mb{G}_2(\tau))
+C\eta\int_{\mathbb{R}}\hspace{-1.15em}-|P_y^{-\mb{X}}|\psi^2dy,
\ea
$$
and
$$
\ba
\int_{\mathbb{R}}\hspace{-1.15em}-
P_y^{-\mb{X}} \f{p_+}{s_3}\f{\phi\z}{v_+\t_+}dy
&=\int_{\mathbb{R}}\hspace{-1.15em}-
P_y^{-\mb{X}} \f{R}{s_3v_+^2}\Big(\phi+\f{\psi}{s_3}-\f{\psi}{s_3}\Big)
\Big(\z-\f{\gamma-1}{R}\f{p_+\psi}{s_3}+\f{\gamma-1}{R}\f{p_+\psi}{s_3}\Big)dy\\
&\leq \int_{\mathbb{R}}\hspace{-1.15em}-
-P_y^{-\mb{X}} (\gamma-1) \f{p_+\psi^2}{s_3^3v_+^2}dy+C\eta\int_{\mathbb{R}}\hspace{-1.15em}-
|P_y^{-\mb{X}}|\psi^2dy +\f{C}{\eta}\sqrt{\d}(\mb{G}_1(\tau)+\mb{G}_2(\tau)).
\ea
$$
We combine the above two estimates, and use (iii) in Lemma \ref{le-shock} and \eqref{change-variable} to have
\begin{align} \label{B1-1}
\mb{B}_{1,1}(\tau) &\leq\int_{\mathbb{R}}-P_y^{-\mb{X}}
\f{\gamma(\gamma+1)}{2}\f{p_+\psi^2}{s_3^3v_+^2}dy
+C(\eta+N(\tau_1,\tau_2)+\sqrt{\d})\int_{\mathbb{R}}|P_y^{-\mb{X}}|\psi^2dy\\ \nm
&\quad+\f{C}{\eta}\sqrt{\d}\big(\mb{G}_1(\tau)+\mb{G}_2(\tau)\big)\\ \nm
&\leq \int_{\mathbb{R}}-P_y^{-\mb{X}}\f{\gamma(\gamma+1)}{2} \f{s_3\psi^2}{\gamma^2p_+} dy +C(\eta+N(\tau_1,\tau_2)+\sqrt{\d})\int_{\mathbb{R}}|P_y^{-\mb{X}}|\psi^2dy\\ \nm
&\quad+\f{C}{\eta}\sqrt{\d}\big(\mb{G}_1(\tau)+\mb{G}_2(\tau)\big)\\ \nm
&=\d\int_0^1\f{\gamma+1}{2\gamma}\f{s_3}{p_+}\psi^2dz
+C(\eta+N(\tau_1,\tau_2)+\sqrt{\d})\d\int_0^1\psi^2dz
+\f{C}{\eta}\sqrt{\d}\big(\mb{G}_1(\tau)+\mb{G}_2(\tau)\big).
\end{align}

Likewise, $\mb{B}_{1,i}(\tau)~(i=2,3)$ can be controlled as
\begin{align} \label{B1-234}
\sum_{i=2}^3\mb{B}_{1,i}(\tau)
&\leq C\d\int_{\mathbb{R}}a^{-\mathbf{X}}_y\psi^2 dy
+C\d\int_{\mathbb{R}}\hspace{-1.15em}-
|P^{-\mathbf{X}}_y|\,|(\phi,\z)|^2dy\\ \nm
&\leq C\d^{\f32}\int_0^1\psi^2dz+C\d^{\f32}\big(\mb{G}_1(\tau)+\mb{G}_2(\tau)\big).
\end{align}
Set
\be\label{alpha+}
\a_+:=\f{\gamma+1}{2\gamma}\f{s_3}{p_+},\quad {\rm with}\quad s_3=\sqrt{-\f{p_--p_+}{v_--v_+}}.
\ee
Hence, combination of \eqref{B1-1} and \eqref{B1-234} gives
\be\label{es-B1}
\mb{B}_1(\tau)\leq \a_+\d\int_0^1\psi^2dz
+C(\eta+N(\tau_1,\tau_2)+\sqrt{\d})\d\int_0^1\psi^2dz
+\f{C}{\eta}\sqrt{\d}\big(\mb{G}_1(\tau)+\mb{G}_2(\tau)\big).
\ee


$\bullet$ {\textbf{Change of variable for $\mathbf{D}_i(\tau)$}}(i=1,2): To deal with the diffusion term $\mb{D}_i(\tau)(i=1,2)$, we need a uniform in $z$ estimate on $\f{dz}{dy}$ with the help of Lemma \ref{le-im}. First, using $a^{-\mathbf{X}}>1$ and the change of variable \eqref{change-variable}, it holds
$$
\ba
\mathbf{D}_1(\tau)+\mathbf{D}_2(\tau)&\geq \int_{\mathbb{R}}\f{\T^{-\mb{X}}}{v\t}\psi_y^2dy
+ \nu\int_{\mathbb{R}}\f{\T^{-\mb{X}}}{v\t^2}\z_y^2dy
\geq\big(1-C(N(\tau_1,\tau_2)+\d)\big)\f{1}{v_+}\int_{\mathbb{R}}\Big(\psi_y^2
+\f{\nu}{\t_+}\z_y^2\Big)dy\\
&=\big(1-C(N(\tau_1,\tau_2)+\d)\big)\f{1}{v_+}\int_0^1\Big(\psi_z^2+\f{\nu}{\t_+}\z_z^2\Big)\f{dz}{dy}dz.
\ea
$$
Integrating \eqref{shock-equ}$_2$ with $\v=1$ over $(-\infty,\xi]$ yields
$$
\f{U^{-\mb{X}}_y}{V^{-\mb{X}}}=s_3^2\big(V^{-\mb{X}}-v_-\big)
+P^{-\mb{X}}-p_-.
$$
On one hand, by \eqref{equivalent} and \eqref{change-variable}, it holds
$$
\f{U^{-\mb{X}}_y}{V^{-\mb{X}}}=\f{1}{s_3}\f{P^{-\mb{X}}_y}{V^{-\mb{X}}}
+\f{ R^s_2}{V^{-\mb{X}}}
=\f{1}{s_3}\f{-\d}{V^{-\mb{X}}}\f{dz}{dy}+\f{ R^s_2}{V^{-\mb{X}}}.
$$
On the other hand,
$$
\ba
&s_3^2\big(V^{-\mb{X}}-v_-\big)+P^{-\mb{X}}-p_-
=-\f{p_+-p_-}{v_+-v_-}\big(V^{-\mb{X}}-v_-\big)+P^{-\mb{X}}-p_-\\
&=\f{1}{v_+-v_-}\Big[\big(P^{-\mb{X}}-p_+\big)\big(V^{-\mb{X}}-v_-\big)
-\big(P^{-\mb{X}}-p_-\big)\big(V^{-\mb{X}}-v_+\big)\Big].
\ea
$$
Recalling $z=(p_- -P^{-\mathbf{X}})/\d$ and $1-z=(P^{-\mathbf{X}}-p_+)/\d$, we have
$$
\f{1}{z(1-z)}\f{dz}{dy}=\f{\d s_3 V^{-\mb{X}}}{ v_+-v_-}\left(\f{V^{-\mathbf{X}}-v_-}{P^{-\mathbf{X}}-p_-}
-\f{V^{-\mathbf{X}}-v_+}{P^{-\mathbf{X}}-p_+}\right)
+\f{ s_3}{z(1-z)}\f{R^s_2}{\d}.
$$
Recalling that
$$
\f{|R^s_2|}{\d}\leq C|P_y^{-\mb{X}}|=C\d\f{dz}{dy},
$$
then using \eqref{important} and \eqref{alpha+}, we obtain
\be\label{coeffi-1}
\left|\f{1}{z(1-z)}\f{dz}{dy}-\d\a_+ v_+\f{R\gamma}{\nu(\gamma-1)^2+R\gamma}\right| \leq C\d^2.
\ee
Thus, by \eqref{coeffi-1} and Poincar\'{e} inequality \eqref{poin-2}, it holds
\be\label{D1D2}
\ba
&\mathbf{D}_1(\tau)+\mathbf{D}_2(\tau)\\
\geq&\big(1-C(N(\tau_1,\tau_2)+\d)\big)
\f{1}{v_+}\int_0^1 z(1-z)\Big(\psi_z^2+\f{\nu}{\t_+}\z_z^2\Big)\f{1}{z(1-z)}\f{dz}{dy}dz\\
\geq& \big(1-C(N(\tau_1,\tau_2)+\d)\big)\d\a_+\f{R\gamma}{\nu(\gamma-1)^2+R\gamma}\int_0^1 z(1-z)\Big(\psi_z^2+\f{\nu}{\t_+}\z_z^2\Big)dz\\
\geq& \big(1-C(N(\tau_1,\tau_2)+\d)\big)\d2\a_+\f{R\gamma}{\nu(\gamma-1)^2+R\gamma}\int_0^1
\Big(\Big|\psi-\int_0^1\psi dz\Big|^2+ \f{\nu}{\t_+}\Big|\z-\int_0^1\z dz\Big|^2\Big) dz\\
=&\big(1-C(N(\tau_1,\tau_2)+\d)\big)\d2\a_+\f{R\gamma}{\nu(\gamma-1)^2+R\gamma}\left(\int_0^1 \psi^2dz-\Big(\int_0^1 \psi dz\Big)^2\right)\\
\quad&+\big(1-C(N(\tau_1,\tau_2)+\d)\big)\d2\a_+\f{R\gamma}{\nu(\gamma-1)^2+R\gamma}\f{\nu}{\t_+}\left(\int_0^1 \z^2dz-\Big(\int_0^1 \z dz\Big)^2\right).
\ea
\ee
Since
$$
1>\f{R\gamma}{\nu(\gamma-1)^2+R\gamma}\rightarrow 0, \quad {\rm as}\quad
\gamma\rightarrow\infty,
$$
we will use the last term in \eqref{D1D2} to control the bad term $\mb{B}_1$. By Cauchy's inequality, it holds
$$
\ba
\int_0^1\z^2dz&=\int_0^1\Big(\z-\f{\gamma-1}{R}\f{p_+\psi}{s_3}+\f{\gamma-1}{R}\f{p_+\psi}{s_3}\Big)^2dz\\
&\geq(1-\d^{\f14})\Big(\f{\gamma-1}{R}\f{p_+}{s_3}\Big)^2\int_0^1\psi^2dz
-C\d^{-\f14}\int_0^1\Big(\z-\f{\gamma-1}{R}\f{p_+\psi}{s_3}\Big)^2dz
\ea
$$
and
$$
\ba
\Big(\int_0^1\z\, dz\Big)^2&\leq 2\Big(\int_0^1\Big(\z-\f{\gamma-1}{R}\f{p_+\psi}{s_3}\Big)dz\Big)^2
+2\Big(\int_0^1\f{\gamma-1}{R}\f{p_+\psi}{s_3} dz\Big)^2\\
&\leq 2\int_0^1\Big(\z-\f{\gamma-1}{R}\f{p_+\psi}{s_3}\Big)^2dz
+2\Big(\f{\gamma-1}{R}\f{p_+}{s_3}\Big)^2\Big(\int_0^1\psi\, dz\Big)^2.
\ea
$$
Thus, using \eqref{change-variable}, it holds
$$
\ba
&\int_0^1\psi^2dz+\f{\nu}{\t_+}\int_0^1\z^2dz\\
\geq & \int_0^1\psi^2dz+(1-\d^{\f14})(\gamma-1)^2\f{\nu\t_+}{s_3^2v_+^2}\int_0^1\psi^2dz
-C\d^{-\f14}\int_0^1\Big(\z-\f{\gamma-1}{R}\f{p_+\psi}{s_3}\Big)^2dz\\
\geq & (1-\d^{\f14})\f{\nu(\gamma-1)^2+R\gamma}{R\gamma}\int_0^1\psi^2dz
-C\d\int_0^1\psi^2dz-C\d^{-\f34}\mb{G}_2(\tau),
\ea
$$
where in the last inequality, we have used (iii) in Lemma \ref{le-shock} that
$$
\f{\t_+}{s_3^2v_+^2}=\f{\t_+}{\gamma p_+ v_+}+O(\d)=\f{1}{R\gamma}+O(\d).
$$
Therefore, we can deduce from \eqref{D1D2} that
\be\label{es-D}
\ba
&\mathbf{D}_1(\tau)+\mathbf{D}_2(\tau)\\
\geq&\big(1-C(N(\tau_1,\tau_2)+\d)\big)\d2\a_+\f{R\gamma}{\nu(\gamma-1)^2+R\gamma}
\Big[(1-\d^{\f14})\f{\nu(\gamma-1)^2+R\gamma}{R\gamma}\int_0^1\psi^2dz\\
\quad&-C\d\int_0^1\psi^2dz-C\d^{-\f34}\mb{G}_2(\tau)
-\f{2\nu(\gamma-1)^2+R\gamma}{R\gamma}\Big(\int_0^1\psi dz\Big)^2
\Big]\\
\geq & \big(1-C(N(\tau_1,\tau_2)+\d^{\f14})\big)\d2\a_+\int_0^1\psi^2dz
-C\d^{\f14}\mb{G}_2(\tau)
-\d2\a_+\f{2\nu(\gamma-1)^2+R\gamma}{\nu(\gamma-1)^2+R\gamma}  \Big(\int_0^1\psi dz\Big)^2.
\ea
\ee


$\bullet$ {\textbf{Estimate on $-\f{\d}{2M}|\dot{\mathbf{X}}(\tau)|^2$}:}
We first decompose $-\f{\d}{M}\dot{\mb{X}}(\tau)$ into
$$
-\f{\d}{M}\dot{\mb{X}}(\tau)=
-\int_{\mathbb{R}}a^{-\mathbf{X}}\f{P^{-\mb{X}}\psi}{V^{-\mb{X}} s_3}V^{-\mathbf{X}}_ydy
+\int_{\mathbb{R}}a^{-\mathbf{X}}\psi U^{-\mathbf{X}}_ydy
+\int_{\mathbb{R}}a^{-\mathbf{X}}\f{p_+\psi}{\T^{-\mb{X}} s_3}\T_y^{-\mathbf{X}}dy\\
=:\sum_{i=4}^6 \mathbf{Y}_{i}(\tau).
$$
Using \eqref{equivalent}, we have
\begin{align}\nm
\mb{Y}_{4}(\tau)&=\int_{\mathbb{R}}a^{-\mathbf{X}}\f{P^{-\mb{X}}\psi}{V^{-\mb{X}} s_3^3}P^{-\mb{X}}_y dy
-\int_{\mathbb{R}}a^{-\mathbf{X}}\f{P^{-\mb{X}}\psi}{V^{-\mb{X}} s_3}R^s_1 dy,\\ \nm
\mb{Y}_{5}(\tau)&=\int_{\mathbb{R}}a^{-\mathbf{X}}\f{\psi}{s_3}P^{-\mb{X}}_y dy
+\int_{\mathbb{R}}a^{-\mathbf{X}}\psi R^s_2 dy,\\ \nm
\mb{Y}_{6}(\tau)&=\int_{\mathbb{R}}a^{-\mathbf{X}}
\f{\gamma-1}{R}\f{p_+\psi}{\T^{-\mb{X}} s_3^3}
P^{-\mb{X}}P^{-\mb{X}}_y dy
+\int_{\mathbb{R}}a^{-\mathbf{X}}\f{p_+\psi}{\T^{-\mb{X}} s_3}R^s_3 dy,
\end{align}
which together with (iii) of Lemma \ref{le-shock} implies
\begin{align}\nm
&\left|\mb{Y}_{4}(\tau)-\int_{\mathbb{R}}\f{s_3v_+}{\gamma^2 p_+}P_y^{-\mb{X}}\psi dy\right|\leq C\sqrt{\d}\int_{\mathbb{R}}|P_y^{-\mb{X}}|\,|\psi|dy,\\ \nm
&\left|\mb{Y}_{5}(\tau)-\int_{\mathbb{R}}\f{s_3v_+}{\gamma p_+}P_y^{-\mb{X}}\psi dy\right|
\leq C\sqrt{\d}\int_{\mathbb{R}}|P_y^{-\mb{X}}|\,|\psi|dy,\\ \nm
&\left|\mb{Y}_{6}(\tau)-(\gamma-1)\int_{\mathbb{R}}\f{s_3v_+}{\gamma^2 p_+}P_y^{-\mb{X}}\psi dy\right|
\leq C\sqrt{\d}\int_{\mathbb{R}}|P_y^{-\mb{X}}|\,|\psi|dy.
\end{align}
Hence, it holds
$$
\left|-\f{\d}{M}\dot{\mb{X}}(\tau)-\int_{\mathbb{R}}\f{2s_3v_+}{\gamma p_+}P_y^{-\mb{X}}\psi dy\right|
\leq C\sqrt{\d}\int_{\mathbb{R}}|P_y^{-\mb{X}}|\,|\psi|dy,
$$
which together with the new change of variable \eqref{change-variable} and \eqref{alpha+} gives
\be\label{es-Y1}
\left|-\f{\d}{M}\dot{\mb{X}}(\tau)+\d\f{4v_+\a_+}{\gamma+1}\int_0^1\psi dz\right|\leq C\d^{\f32}\int_0^1|\psi|dz.
\ee
By \eqref{es-Y1}, it holds
$$
\ba
\left|\dot{\mathbf{X}}(\tau)-M\f{4v_+\a_+}{\gamma+1}\int_0^1 \psi dz\right|
=\left|\f{M}{\d}\Big(-\f{\d}{M}\dot{\mb{X}}(\tau)+\d\f{4v_+\a_+}{\gamma+1}\int_0^1\psi dz\Big)\right|
\leq C\sqrt{\d}\int_0^1|\psi|dz,
\ea
$$
which yields
$$
\ba
\left(\left|M\f{4v_+\a_+}{\gamma+1}\int_0^1 \psi dz\right|
-|\dot {\mathbf{X}}(\tau)|\right)^2
&\leq C\d\int_0^1 \psi^2 dz,
\ea
$$
which together with the algebraic inequality $\f{m^2}{2}-n^2\leq (m-n)^2$  indicates
$$
M^2\f{8v_+^2\alpha_+^2}{(\gamma+1)^2}\left(\int_0^1 \psi dz\right)^2
-|\dot{\mathbf{X}}(\tau)|^2\leq C\d\int_0^1 \psi^2 dz.
$$
Thus, we get
\be\label{es-X}
-\f{\d}{2M}|\dot{\mathbf{X}}(\tau)|^2\leq
-M\f{4v_+^2\a_+^2}{(\gamma+1)^2}\d\left(\int_0^1 \psi dz\right)^2
+C\d^2\int_0^1 \psi^2 dz.
\ee


$\bullet$ \textbf{Conclusion}: First, combination of \eqref{es-B1}, \eqref{es-D} and \eqref{es-X} yields
$$
\ba
&\quad\mathbf{B}_1(\tau)-\f78\big(\mathbf{D}_1(\tau)+\mathbf{D}_2(\tau)\big)-\f{\d}{2M}|\dot{\mb{X}}(\tau)|^2\\
&\leq \a_+\d\int_0^1\psi^2dz
+C(\eta+N(\tau_1,\tau_2)+\sqrt{\d})\d\int_0^1\psi^2dz
+\f{C}{\eta}\sqrt{\d}\big(\mb{G}_1(\tau)+\mb{G}_2(\tau)\big)\\
&-\big(1-C(N(\tau_1,\tau_2)+\d^\f{1}{4})\big)\f74\a_+\d\int_0^1 \psi^2dz+C\d^{\f14}\mb{G}_2(\tau)\\
&+2\a_+\f{2\nu(\gamma-1)^2+R\gamma}{\nu(\gamma-1)^2+R\gamma}
\d\Big(\int_0^1\psi dz\Big)^2
-M\f{4v_+^2\a_+^2}{(\gamma+1)^2}\d\Big(\int_0^1 \psi dz\Big)^2.
\ea
$$
Choosing $N(\tau_1,\tau_2)$, $\d$ and $\eta$ suitably small such that
$$
C(\eta+N(\tau_1,\tau_2)+\sqrt{\d})\leq \f{\alpha_+}{3},\quad
C(N(\tau_1,\tau_2)+\d^{\f14})\leq \f{1}{21},
$$
which together with $M=\f{(\gamma+1)^2}{2\a_+v_+^2}\f{2\nu(\gamma-1)^2+R\gamma}{\nu(\gamma-1)^2+R\gamma}$ gives
$$
\mathbf{B}_1(\tau)-\f78\big(\mathbf{D}_1(\tau)+\mathbf{D}_2(\tau)\big)-\f{\d}{2M}|\dot{\mb{X}}(\tau)|^2
\leq -\f{\a_+}{3}\d\int_0^1\psi^2dz+C\d^{\f14}\big(\mb{G}_1(\tau)+\mb{G}_2(\tau)\big),
$$
which indicates the desired inequality \eqref{es-contaction} by using \eqref{change-variable}. The proof of Lemma \ref{le-contraction} is completed.

\hfill $\Box$

%
%
%

\subsection{Energy estimates}

Now we are in a stage to show the $L^2$ estimates.

\begin{lemma}\label{le-L2}
Under the hypotheses of Proposition \ref{priori}, there exists constant $C>0$,
such that for any $\tau\in[\tau_1,\tau_2]$, it holds
\begin{equation}\label{basic}
\ba
&\|(\phi,\psi,\z)(\tau)\|^2+\int_{\tau_1}^\tau
\Big(\d|\dot{\mathbf{X}}(\tau)|^2+G^{s}(\tau)+\sum_{i=1}^2G_{i}(\tau)\Big)d\tau
+\int_{\tau_1}^\tau\|(\psi_y,\z_y)\|^2d\tau\\
&\leq C\|(\phi,\psi,\z)(\tau_1)\|^2,
\ea
\end{equation}
where we denote
\be\label{new-G...}
\ba
G^{s}(\tau)&=\int_{\mathbb{R}}|U_y^{-\mathbf{X}}|\psi^2dy,\quad
G_{1}(\tau)=\f{1}{\sqrt{\d}}\int_{\mathbb{R}}\hspace{-1.15em}-
|U_y^{-\mathbf{X}}|\left(\phi+\f{\psi}{s_3}\right)^2dy,\\
G_{2}(\tau)&=\f{1}{\sqrt{\d}}\int_{\mathbb{R}}\hspace{-1.15em}-
|U_y^{-\mathbf{X}}|\left(\z-\f{\gamma-1}{R}\f{p_+\psi}{s_3}\right)^2dy.
\ea
\ee

\end{lemma}

\textbf{\emph{Proof}}:  First of all, we use \eqref{ba-expression} and Lemma \ref{le-contraction} to have
\be\label{pre-basic}
\ba
&\quad\f{d}{d\tau}\int_{\mathbb{R}}\hspace{-1.15em}-
a^{-\mathbf{X}}\mathcal{E}dy+\mathbf{G}_1(\tau)+\mathbf{G}_2(\tau)
+\f18\big(\mathbf{D}_1(\tau)+\mathbf{D}_2(\tau)\big)+\f{\d}{2M}|\dot{\mathbf{X}}(\tau)|^2\\
&=\mb{B}_1(\tau)-\f78\big(\mb{D}_1(\tau)+\mb{D}_2(\tau)\big)-\f{\d}{2M}|\dot{\mathbf{X}}(\tau)|^2
+\sum_{i=2}^{7}\mathbf{B}_{i}(\tau)
+\dot{\mathbf{X}}(\tau)\sum_{i=1}^3\mathbf{Y}_i(\tau)\\
&\leq -\f{\alpha_+}{3}\underbrace{\int_{\mathbb{R}}|P^{-\mathbf{X}}_y|\psi^2dy}_{\mb{G}^s(\tau)}
+C\d^{\f14}\big(\mathbf{G}_1(\tau)+\mathbf{G}_2(\tau)\big)
+\sum_{i=2}^{7}\mathbf{B}_{i}(\tau)+\f{\d}{4M}|\dot{\mathbf{X}}(\tau)|^2
+\f{C}{\d}\sum_{i=1}^3|\mathbf{Y}_i(\tau)|^2.
\ea
\ee
In what follows, it is to control the bad terms $\mathbf{B}_{i}(\tau)\,(i=2,\cdots,7)$ and $\mathbf{Y}_{i}(\tau)\,(i=1,2,3)$ by the good terms $\mathbf{G}^{s}(\tau)$, $\mathbf{G}_{i}(\tau)$ and $\mathbf{D}_i(\tau)~(i=1,2)$.

\

$\bullet$ \textbf{Estimates for $\mathbf{B}_{i}(\tau)~(i=2,\cdots,7)$}: The smallness of assumption \eqref{assumption} and wave strengths leads to
$$
\mathbf{B}_{2}(\tau)\leq C(N(\tau_1,\tau_2)+\d)\big(\mb{G}_1(\tau)+\mb{G}_2(\tau)\big).
$$

For $\mathbf{B}_{3}(\tau)$, it holds
\begin{align}\nm
\mb{B}_3(\tau)&\leq C \int_{\mathbb{R}}\hspace{-1.15em}-a^{-\mathbf{X}}_y
\left|\left(\phi+\f{\psi}{s_3}, \z-\f{\gamma-1}{R}\f{p_+\psi}{s_3}\right)\right|\psi^2dy
+\f{C}{\sqrt{\d}}\int_{\mathbb{R}}|P^{-\mb{X}}_y|\,|\psi|^3dy\\ \nm
&+C\sqrt{\d}\int_{\mathbb{R}}|P^{-\mb{X}}_y|\psi^2dy
=:\sum_{i=1}^2\mb{B}_{3,i}(\tau)+C\sqrt{\d}\mb{G}^s(\tau).
\end{align}
By Cauchy inequality, Gagliardo-Nirenberg inequality and \eqref{a>1}, it holds
$$
\ba
\mb{B}_{3,1}(\tau)&\leq C\|\psi\|^2_{L^{\infty}}\sqrt{\mathbf{G}_{1}(\tau)+\mathbf{G}_{2}(\tau)}
\Big(\int_{\mathbb{R}}|a_y^{-\mathbf{X}}|dy\Big)^{\f12}\\
&\leq C\|\psi\|\|\psi_y\|\d^{\f14}\sqrt{\mathbf{G}_{1}(\tau)+\mathbf{G}_{2}(\tau)}
\leq C\d^{\f14}N(\tau_1,\tau_2)\big(\|\psi_y\|^2+\mathbf{G}_{1}(\tau)+\mathbf{G}_{2}(\tau)\big).
\ea
$$
Likewise, we have
$$
\ba
\mb{B}_{3,2}(\tau)&\leq \f{C}{\sqrt{\d}}\|\psi\|^2_{L^{\infty}}\int_{\mathbb{R}}|P_y^{-\mathbf{X}}|\,|\psi|dy
\leq \f{C}{\sqrt{\d}}\|\psi\|\|\psi_y\|
\sqrt{\mathbf{G}^{s}(\tau)}\Big(\int_{\mathbb{R}}|P_y^{-\mathbf{X}}|dy\Big)^{\f12}\\
&\leq CN(\tau_1,\tau_2)\|\psi_y\|\sqrt{\mathbf{G}^{s}(\tau)}
\leq CN(\tau_1,\tau_2)\big(\|\psi_y\|^2+\mathbf{G}^{s}(\tau)\big).
\ea
$$
Combination the estimates of $\mb{B}_{3,i}(\tau)~(i=1,2)$ gives
$$
|\mathbf{B}_{3}(\tau)|\leq  CN(\tau_1,\tau_2)\big(\|\psi_y\|^2+\mathbf{G}_{1}(\tau)+\mathbf{G}_{2}(\tau)+\mathbf{G}^{s}(\tau)\big).
$$
For $\mathbf{B}_{4}(\tau)$ and $\mathbf{B}_{5}(\tau)$, using Cauchy inequality, we have
$$
\ba
|\mathbf{B}_{4}(\tau)|+|\mathbf{B}_{5}(\tau)|&\leq \f{1}{40}\int_{\mathbb{R}}a^{-\mb{X}}\f{\T^{-\mb{X}}\psi_y^2}{v\t}dy
+\f{\nu}{40}\int_{\mathbb{R}}a^{-\mb{X}}\f{\T^{-\mb{X}}\z_y^2}{v\t^2}dy
+C\int_{\mathbb{R}}\hspace{-1.15em}-
|(U^{-\mb{X}}_y,\T^{-\mb{X}}_y)|^2|(\phi,\z)|^2dy\\
&\leq \f{1}{40}\big(\mb{D}_1(\tau)+\mb{D}_2(\tau)\big)
+C\d^2\mb{G}^s(\tau)+C\d^{\f52}\big(\mathbf{G}_{1}(\tau)+\mathbf{G}_{2}(\tau)\big).
\ea
$$
Likewise,
$$
\ba
|\mathbf{B}_{6}(\tau)|&\leq \f{1}{40}\big(\mb{D}_1(\tau)+\mb{D}_2(\tau)\big)
+C\int_{\mathbb{R}}|a^{-\mb{X}}_y|^2|(\psi,\z)|^2dy,
\ea
$$
where
$$
\ba
C\int_{\mathbb{R}}|a^{-\mb{X}}_y|^2|(\psi,\z)|^2dy
&\leq C\|a^{-\mb{X}}_y\|_{L^{\infty}}\mb{G}_2(\tau)
+\f{C}{\d}\|P^{-\mb{X}}_y\|_{L^{\infty}}\mb{G}^s(\tau)\\
&\leq C\d^{\f32}\mb{G}_2(\tau)+C\d\mb{G}^s(\tau),
\ea
$$
which indicates
$$
|\mathbf{B}_{6}(\tau)|\leq \f{1}{40}\big(\mb{D}_1(\tau)+\mb{D}_2(\tau)\big)
+C\d^{\f32}\mb{G}_2(\tau)+C\d\mb{G}^s(\tau).
$$
By (ii) of Lemma \ref{le-shock} and \eqref{a>1}, it holds
$$
\ba
\mathbf{B}_{7}(\tau)&\leq C\int_{\mathbb{R}}\hspace{-1.15em}-a_y^{-\mb{X}}|(U^{-\mb{X}}_y,\T^{-\mb{X}}_y)|\,|(\phi,\psi,\z)|^2dy\\
&\leq C\|(U_y^{-\mb{X}},\T_y^{-\mb{X}})\|_{L^{\infty}}
\int_{\mathbb{R}}\hspace{-1.15em}-a_y^{-\mb{X}}\left|\left(\phi+\f{\psi}{s_3}, \z-\f{\gamma-1}{R}\f{p_+\psi}{s_3}\right)\right|^2dy\\
&\quad+C\|a_y^{-\mb{X}}\|_{L^{\infty}}\int_{\mathbb{R}}|(U_y^{-\mb{X}},\T_y^{-\mb{X}})|\psi^2dy\\
&\leq C\d^2\big(\mathbf{G}_{1}(\tau)+\mathbf{G}_{2}(\tau)\big)+C\d^{\f32}\mb{G}^s(\tau).
\ea
$$
Thus, we combine the above estimates and choose $N(\tau_1,\tau_2)$ and $\d$   suitably small such that
\be\label{sum-B}
\ba
\sum_{i=2}^{7}\mathbf{B}_{i}(\tau)
&\leq \f{1}{16}\big(\mathbf{G}^{s}(\tau)+\mathbf{G}_{1}(\tau)
+\mathbf{G}_{2}(\tau)+\mathbf{D}_1(\tau)+\mathbf{D}_2(\tau)\big).
\ea
\ee


$\bullet$ \textbf{Estimates for $\mathbf{Y}_{i}(\tau)~(i=1,2,3)$}: For $\mathbf{Y}_{1}(\tau)$,  it holds
$$
\ba
\f{C}{\d}|\mathbf{Y}_{1}(\tau)|^2&\leq C\left|\int_{\mathbb{R}}\hspace{-1.15em}-
a_y^{-\mathbf{X}}\Big|\Big(\phi+\f{\psi}{s_3}, \z-\f{\gamma-1}{R}\f{p_+\psi}{s_3}\Big)\Big|dy\right|^2\\
&\leq C\big(\mb{G}_1(\tau)+\mb{G}_2(\tau)\big)\int_{\mathbb{R}}a_y^{-\mb{X}}dy
\leq C\sqrt{\d}\big(\mb{G}_1(\tau)+\mb{G}_2(\tau)\big).
\ea
$$
For $\mathbf{Y}_{2}(\tau)$, it holds
$$
\ba
\f{C}{\d}|\mathbf{Y}_{2}(\tau)|^2&\leq \f{C}{\d}\left|
\int_{\mathbb{R}}\hspace{-1.15em}-
|\T_y^{-\mb{X}}| \, |(\phi,\z)|^2dy\right|^2\\
&\leq \f{C}{\d}\|\T_y^{-\mb{X}}\|_{L^{\infty}} \|(\phi,\z)\|^2
\int_{\mathbb{R}}\hspace{-1.15em}-|\T_y^{-\mb{X}}| \, |(\phi,\z)|^2dy\\
&\leq C\d N(\tau_1,\tau_2)^2\left(\sqrt{\d}\big(\mb{G}_1(\tau)+\mb{G}_2(\tau)\big)
+\mb{G}^s(\tau)\right).
\ea
$$
Likewise, we have
$$
\ba
\f{C}{\d}|\mathbf{Y}_{3}(\tau)|^2&\leq \f{C}{\d}\left|
\int_{\mathbb{R}}\hspace{-1.15em}-a_y^{-\mb{X}} \, |(\phi,\psi,\z)|^2dy\right|^2
\leq \f{C}{\d}\|a_y^{-\mb{X}}\|_{L^{\infty}} \|(\phi,\psi,\z)\|^2
\int_{\mathbb{R}}\hspace{-1.15em}-a_y^{-\mb{X}} \, |(\phi,\psi,\z)|^2dy\\
&\leq C\sqrt{\d}N(\tau_1,\tau_2)^2\left(\mb{G}_1(\tau)+\mb{G}_2(\tau)
+\f{1}{\sqrt{\d}}\mb{G}^s(\tau)\right).
\ea
$$
Thus, we combine the above estimates and choose $N(\tau_1,\tau_2)$ and $\d$ suitably small to have
\be\label{sum-Y}
\sum_{i=1}^3\f{C}{\d}|\mathbf{Y}_{i}(\tau)|^2\leq \f{1}{16}
\big(\mathbf{G}_{1}(\tau)+\mathbf{G}_{2}(\tau)+\mathbf{G}^{s}(\tau)\big).
\ee
We substitute \eqref{sum-B} and \eqref{sum-Y} into \eqref{pre-basic}, then integrate the resultant inequality over $(\tau_1,\tau)$, and note that
$$
\mb{G}^s(\tau)\sim G^s(\tau), \quad \mb{G}_i(\tau)\sim G_i(\tau),
\quad \mb{D}_1(\tau)\sim \|\psi_y\|^2,\quad \mb{D}_2(y)\sim \|\z_y\|^2,
$$
we can obtain the desired inequality \eqref{basic} with the new notations \eqref{new-G...}. Here $A \sim B$ denotes $A$ is equivalent to $B$. The proof of Lemma \ref{le-L2} is completed.

\hfill $\Box$

%
%
%
%

Next, we estimate the derivative of $\phi_y$.

\begin{lemma}\label{le-1st}
Under the hypotheses of Proposition \ref{priori}, there exists constant $C>0$,
such that for any $\tau\in[\tau_1,\tau_2]$, it holds
\be\label{v-1st}
\ba
\|\phi_y(\tau)\rlap{--}\hspace{0.07em}\|^2+\int_{\tau_1}^\tau\|\phi_y\rlap{--}\hspace{0.07em}\|^2d\tau\leq C(\|(\phi,\psi,\z)(\tau_1)\|^2+\|\phi_{y}(\tau_1)\rlap{--}\hspace{0.07em}\|^2+\d).
\ea
\ee
\end{lemma}

\textbf{\emph{Proof}}: Set $\tilde v=\f{v}{V^{-\mb{X}}}$, then the equation \eqref{perturb}$_2$ can be rewritten as
\be\label{phi-y}
\ba
&\Big(\f{\tilde v_y}{\tilde v}\Big)_{\tau}-\dot{\mb{X}}(\tau)\left(\f{V^{-\mb{X}}_{yy}}{V^{-\mb{X}}}
-\f{(V^{-\mb{X}}_{y})^2}{(V^{-\mb{X}})^2}-U^{-\mb{X}}_y\right)-\psi_{\tau}
+\f{R\t}{v}\f{\tilde v_y}{\tilde v}\\
&-\f{R\z_y}{v}+(p-P^{-\mb{X}})\f{V^{-\mb{X}}_y}{V^{-\mb{X}}}
-R\T^{-\mb{X}}_y\left(\f{1}{v}-\f{1}{V^{-\mb{X}}}\right)=0,
\ea
\ee
where we have used the fact
$$
-\big(p-P^{-\mb{X}}\big)_y=\f{R\t}{v}\f{\tilde v_y}{\tilde v}
-\f{R\z_y}{v}+(p-P^{-\mb{X}})\f{V^{-\mb{X}}_y}{V^{-\mb{X}}}
-R\T^{-\mb{X}}_y\left(\f{1}{v}-\f{1}{V^{-\mb{X}}}\right).
$$
Multiplying \eqref{phi-y} by $\f{\tilde v_y}{\tilde v}$, we have
$$
\ba
&\left(\f12\Big(\f{\tilde v_y}{\tilde v}\Big)^2-\psi\f{\tilde v_y}{\tilde v}\right)_{\tau}+\left(\psi\f{\tilde v_\tau}{\tilde v}\right)_{y}
+\f{R\t}{v}\Big(\f{\tilde v_y}{\tilde v}\Big)^2\\
=&\dot{\mb{X}}(\tau)\left(\f{V^{-\mb{X}}_{yy}}{V^{-\mb{X}}}
-\f{(V^{-\mb{X}}_{y})^2}{(V^{-\mb{X}})^2}-U^{-\mb{X}}_y\right)\f{\tilde v_y}{\tilde v}+\dot{\mb{X}}(\tau)\f{V^{-\mb{X}}_y}{V^{-\mb{X}}}\psi_y
+\psi_y\left(\f{u_y}{v}-\f{U^{-\mb{X}}_y}{V^{-\mb{X}}}\right)\\
+&\left(\f{R\z_y}{v}-(p-P^{-\mb{X}})\f{V^{-\mb{X}}_y}{V^{-\mb{X}}}
+R\T^{-\mb{X}}_y\Big(\f{1}{v}-\f{1}{V^{-\mb{X}}}\Big)
\right)\f{\tilde v_y}{\tilde v}.
\ea
$$
Integrating the above equality with respect to $\tau$ and $y$ over $[\tau_1,\tau]\times\mathbb{R}_{\pm}$ and using Cauchy inequality, we obtain
\be\label{int-phiy}
\ba
&\int_{\mathbb{R}}\hspace{-1.15em}-\left(\f12\Big(\f{\tilde v_y}{\tilde v}\Big)^2-\psi\f{\tilde v_y}{\tilde v}\right)(\tau,y)dy
+\int_{\tau_1}^{\tau}\int_{\mathbb{R}}\hspace{-1.15em}-\f{R\t}{2v}\Big(\f{\tilde v_y}{\tilde v}\Big)^2dyd\tau\\
\leq& \int_{\mathbb{R}}\hspace{-1.15em}-\left(\f12\Big(\f{\tilde v_y}{\tilde v}\Big)^2-\psi\f{\tilde v_y}{\tilde v}\right)(\tau_1,y)dy
+\int_{\tau_1}^{\tau}\left[\psi\f{\tilde v_{\tau}}{\tilde v}\right](\tau)d\tau
+\int_{\tau_1}^{\tau}\int_{\mathbb{R}}\hspace{-1.15em}-
\left|\psi_y\Big(\f{u_y}{v}-\f{U^{-\mb{X}}}{V^{-\mb{X}}}\Big)\right|dyd\tau\\
+&\int_{\tau_1}^{\tau}|\dot{\mb{X}}(\tau)|\int_{\mathbb{R}}
\left|\f{V^{-\mb{X}}_y}{V^{-\mb{X}}}\psi_y\right|dyd\tau
+C\int_{\tau_1}^{\tau}|\dot{\mb{X}}(\tau)|^2\int_{\mathbb{R}}
\left(\f{V^{-\mb{X}}_{yy}}{V^{-\mb{X}}}
-\f{(V^{-\mb{X}}_{y})^2}{(V^{-\mb{X}})^2}-U^{-\mb{X}}_y\right)^2dyd\tau\\
+&C\int_{\tau_1}^{\tau}\int_{\mathbb{R}}\hspace{-1.15em}-
\left|\f{R\z_y}{v}-(p-P^{-\mb{X}})\f{V^{-\mb{X}}_y}{V^{-\mb{X}}}
+R\T^{-\mb{X}}_y\Big(\f{1}{v}-\f{1}{V^{-\mb{X}}}\Big)\right|^2dyd\tau.
\ea
\ee
To deal with the terms in \eqref{int-phiy} which involve jumps at $y=0$, we recall the observation by Hoff \cite{Hoff} that the jump discontinuity of the solution $(v,u)$ to \eqref{new-euler} and \eqref{new-data} satisfies
$$
[\log v(\tau,0)]=[\log v(\tau_1,0)]\exp\Big\{\int_{\tau_1}^{\tau}\alpha(s)ds\Big\}
$$
with
$$
\alpha(\tau)=\f{p\big(\log v(\tau,0+),\t(\tau,0)\big)
-p\big(\log v(\tau,0-),\t(\tau,0)\big)}{\log v(\tau,0+)-\log v(\tau,0-)},
$$
from which it follows for $\tau\geq\tau_1$,
\be\label{es-jump}
\left[\f{u_y}{v}\right]=[p],\quad |\,[v](\tau)\,|\leq C|\,[v](\tau_1)\,|\exp\{-C(\tau-\tau_1)\}.
\ee
Hence, by \eqref{es-jump}, the jump across $y=0$ can be bounded as follows
$$
\ba
&\int_{\tau_1}^{\tau}\left[\psi\f{\tilde v_\tau}{\tilde v}\right](\tau)d\tau
=\int_{\tau_1}^{\tau}\psi(\tau,0)\left[\f{v_\tau}{v}-\f{V_\tau^{-\mb{X}}}{V^{-\mb{X}}}\right](\tau)d\tau\\
=&\int_{\tau_1}^{\tau}\psi(\tau,0)\left[\f{u_y}{v}-\f{U_y^{-\mb{X}}}{V^{-\mb{X}}}\right](\tau)d\tau
+\int_{\tau_1}^{\tau}\psi(\tau,0)\dot{\mb{X}}(\tau)\left[\f{V_y^{-\mb{X}}}{V^{-\mb{X}}}\right](\tau)d\tau\\
=&\int_{\tau_1}^{\tau}\psi(\tau,0)[p](\tau)d\tau
=R\int_{\tau_1}^{\tau}\psi(\tau,0)\t(\tau,0)\left[\f{1}{v}\right](\tau)d\tau
=-R\int_{\tau_1}^{\tau}\f{\psi(\tau,0)\t(\tau,0)}{v(\tau,0+)v(\tau,0-)}[v](\tau)d\tau\\
\leq& C\int_{\tau_1}^{\tau}\|\psi\|_{L^{\infty}}(\tau)|[v](\tau_1)|e^{-C(\tau-\tau_1)}d\tau
\leq C\d\int_{\tau_1}^{\tau}\|\psi\|^{\f12}\|\psi_y\|^{\f12}e^{-C(\tau-\tau_1)}d\tau\\
\leq & C\d \int_{\tau_1}^{\tau}\|\psi_y\|^2d\tau
+C\d\sup_{\tau\in[\tau_1,\tau_2]}\|\psi\|^2(\tau)+C\d.
\ea
$$
By Cauchy inequality, it holds
$$
\ba
\int_{\tau_1}^{\tau}\int_{\mathbb{R}}\hspace{-1.15em}-
\left|\psi_y\Big(\f{u_y}{v}-\f{U^{-\mb{X}}_y}{V^{-\mb{X}}}\Big)\right|dyd\tau
&\leq C\int_{\tau_1}^{\tau}\|\psi_y\|^2d\tau+C\int_{\tau_1}^{\tau}
\int_{\mathbb{R}}\hspace{-1.15em}-|U^{-\mb{X}}_y|^2\phi^2dyd\tau\\
&\leq  C\int_{\tau_1}^{\tau}\|\psi_y\|^2d\tau
+C\d^2\int_{\tau_1}^{\tau}\big(G^s(\tau)+G_1(\tau)\big)d\tau,
\ea
$$
and
$$
\ba
&\int_{\tau_1}^{\tau}|\dot{\mb{X}}(\tau)|\int_{\mathbb{R}}
\left|\f{V^{-\mb{X}}_y}{V^{-\mb{X}}}\psi_y\right|dyd\tau
+C\int_{\tau_1}^{\tau}|\dot{\mb{X}}(\tau)|^2\int_{\mathbb{R}}
\left(\f{V^{-\mb{X}}_{yy}}{V^{-\mb{X}}}
-\f{(V^{-\mb{X}}_{y})^2}{(V^{-\mb{X}})^2}-U^{-\mb{X}}_y\right)^2dyd\tau\\
&\leq C\int_{\tau_1}^{\tau}|\dot{\mb{X}}(\tau)|\,\|V^{-\mb{X}}_y\|\,\|\psi_y\|d\tau
+C\int_{\tau_1}^{\tau}|\dot{\mb{X}}(\tau)|^2\,\|(V^{-\mb{X}}_{yy},(V^{-\mb{X}}_y)^2,U^{-\mb{X}}_y)\|^2d\tau\\
&\leq C\int_{\tau_1}^{\tau}\|\psi_y\|^2d\tau
+C\d^3\int_{\tau_1}^{\tau}|\dot{\mb{X}}(\tau)|^2d\tau.
\ea
$$
Similarly, we have
$$
\ba
&C\int_{\tau_1}^{\tau}\int_{\mathbb{R}}\hspace{-1.15em}-
\left|\f{R\z_y}{v}-(p-P^{-\mb{X}})\f{V^{-\mb{X}}_y}{V^{-\mb{X}}}
-R\T^{-\mb{X}}_y\Big(\f{1}{v}-\f{1}{V^{-\mb{X}}}\Big)\right|^2dyd\tau\\
&\leq C\int_{\tau_1}^{\tau}\|\z_y\|^2d\tau
+C\int_{\tau_1}^{\tau}\int_{\mathbb{R}}\hspace{-1.15em}-
|(V^{-\mb{X}}_y,\T^{-\mb{X}}_y)|^2|(\phi,\z)|^2dyd\tau\\
&\leq C\int_{\tau_1}^{\tau}\|\z_y\|^2d\tau
+C\d^2\int_{\tau_1}^{\tau}\big(G^s(\tau)+G_1(\tau)+G_2(\tau)\big)d\tau.
\ea
$$
Note that
$$
\f{\tilde v_y}{\tilde v}=\f{v_y}{v}-\f{V^{-\mb{X}}_y}{V^{-\mb{X}}}
=\f{\phi_y}{v}-\f{\phi V^{-\mb{X}}_y}{v V^{-\mb{X}}},
$$
which yields
$$
C^{-1}\big(\phi_y^2-|V^{-\mb{X}}_y\phi|^2\big)\leq \Big(\f{\tilde v_y}{\tilde v}\Big)^2\leq C\big(\phi_y^2+|V^{-\mb{X}}_y\phi|^2\big).
$$
Thus, combination of the above inequalities and Lemma \ref{le-L2} implies the desired inequality \eqref{v-1st}. The proof of Lemma \ref{le-1st} is completed.

\hfill $\Box$

\begin{lemma}\label{le-1st-ut}
Under the hypotheses of Proposition \ref{priori}, there exists constant $C>0$,
such that for any $\tau\in[\tau_1,\tau_2]$, it holds
\be\label{ut-1st}
\ba
\|(\psi_y,\z_y)(\tau)\|^2+\int_{\tau_1}^\tau(\|(\psi_{yy},\z_{yy})\rlap{--}\hspace{0.07em}\|^2
+\|(\psi_{\tau},\z_{\tau})\|^2)d\tau
\leq C(N^2(\tau_1)+\d)+C\d\int_{\tau_1}^{\tau}\|(\psi_{\tau y},\z_{\tau y})\rlap{--}\hspace{0.07em}\|^2d\tau.
\ea
\ee
\end{lemma}

\textbf{\emph{Proof}}: 
We multiply \eqref{perturb}$_2$ by $-\psi_{yy}$ to have
\be\label{u-1st}
\ba
&\left(\f{\psi_y^2}{2}\right)_\tau-(\psi_y\psi_\tau)_y+\f{\psi_{yy}^2}{v}\\
&=-\dot{\mathbf{X}}(\tau)U^{-\mathbf{X}}_y\psi_{yy}
+\left(-\f{p}{v}\phi_y+\f{R}{v}\z_y \right)\psi_{yy}
+R\T^{-\mb{X}}_y\left(\f{1}{v}-\f{1}{V^{-\mb{X}}}\right)\psi_{yy}\\
&-RV^{-\mb{X}}_y\left(\f{\t}{v^2}-\f{\T^{-\mb{X}}}{(V^{-\mb{X}})^2}\right)\psi_{yy}
-U^{-\mb{X}}_{yy}\left(\f{1}{v}-\f{1}{V^{-\mb{X}}}\right)\psi_{yy}
+\f{1}{v^2}\phi_y\psi_y\psi_{yy}\\
&+\f{1}{v^2}U^{-\mb{X}}_y\phi_y\psi_{yy}+\f{1}{v^2}V^{-\mb{X}}_y\psi_y\psi_{yy}+V^{-\mb{X}}_y U^{-\mb{X}}_y\left(\f{1}{v^2}-\f{1}{(V^{-\mb{X}})^2}\right)\psi_{yy}.
\ea
\ee
Integrating \eqref{u-1st} with respect to $\tau$ and $y$ over $[\tau_1,\tau]\times\mathbb{R}_{\pm}$ yields
\be\label{u-x}
\ba
&\int_{\mathbb{R}}\f{\psi_y^2}{2}(\tau,y)dy
+\int_{\tau_1}^{\tau}\int_{\mathbb{R}}\hspace{-1.15em}-\f{\psi_{yy}^2}{2v}dyd\tau
\leq \int_{\mathbb{R}}\f{\psi_y^2}{2}(\tau_1,y)dy
-\int_{\tau_1}^{\tau}[\psi_y\psi_{\tau}](\tau)d\tau\\
&+C\int_{\tau_1}^{\tau}|\dot{\mb{X}}(\tau)|^2\int_{\mathbb{R}}|U^{-\mb{X}}_{y}|^2dyd\tau
+C\int_{\tau_1}^{\tau}\int_{\mathbb{R}}\hspace{-1.15em}-(\phi_y^2+\z_y^2)dyd\tau
+C\int_{\tau_1}^{\tau}\int_{\mathbb{R}}\hspace{-1.15em}-|\phi_y|\,|\psi_y|\,|\psi_{yy}|dyd\tau\\
&+C\int_{\tau_1}^{\tau}\int_{\mathbb{R}}\hspace{-1.15em}-|(U^{-\mb{X}}_{yy},V^{-\mb{X}}_y,\T^{-\mb{X}}_y)|^2|(\phi,\z)|^2dyd\tau
+C\int_{\tau_1}^{\tau}\int_{\mathbb{R}}\hspace{-1.15em}-|(V^{-\mb{X}}_y,U^{-\mb{X}}_y)|^2|(\phi_y,\psi_y)|^2dyd\tau.
\ea
\ee
We control the jump term across $y=0$ as follows,
\be\label{jump-3}
\ba
-\int_{\tau_1}^{\tau}[\psi_y\psi_{\tau}](\tau)d\tau
&=-\int_{\tau_1}^{\tau}\psi_{\tau}(\tau,0)[\psi_y](\tau)d\tau
=-\int_{\tau_1}^{\tau}\psi_{\tau}(\tau,0)[u_y](\tau)d\tau\\
&=-\int_{\tau_1}^{\tau}\psi_{\tau}(\tau,0)\left[\Big(\f{u_y}{v}-p\Big)v\right](\tau)d\tau\\
&=-\int_{\tau_1}^{\tau}\psi_{\tau}(\tau,0)\Big(\f{u_y}{v}-p\Big)(\tau,0)[v](\tau)d\tau\\
&\leq C\int_{\tau_1}^{\tau}\|\psi_{\tau}\|_{L^{\infty}}
(\|\psi_{y}\|_{L^{\infty}}+1)[v](\tau_1)e^{-c(\tau-\tau_1)}d\tau\\
&\leq C\d\int_{\tau_1}^{\tau}\|\psi_{\tau}\|^{\f12}\|\psi_{\tau y}\rlap{--}\|^{\f12}
(\|\psi_{y}\|^{\f12}\|\psi_{yy}\rlap{--}\|^{\f12}+1)e^{-c(\tau-\tau_1)}d\tau.
\ea
\ee
By equation \eqref{perturb}$_2$ and assumption \eqref{assumption}, it holds
\be\label{psi_tau}
\ba
\|\psi_{\tau}\|&\leq C\big(|\dot{\mb{X}}(\tau)|\|U^{-\mb{X}}_y\|+\|\psi_{yy}\rlap{--}\|+\|(\phi_y,\psi_y,\z_y)\rlap{--}\|
+\|\phi_y\psi_y\rlap{--}\|\\
&\quad+\|(V^{-\mb{X}}_y,U^{-\mb{X}}_y,\T^{-\mb{X}}_y,U^{-\mb{X}}_{yy})\phi\|+\|V^{-\mb{X}}_y\z\|\big)
&\\
&\leq C\big(\d^{\f32}|\dot{\mb{X}}(\tau)|+\|\psi_{yy}\rlap{--}\|+\|(\phi_y,\psi_y,\z_y)\rlap{--}\|
+\d^{\f32}\big),
\ea
\ee
where we have used the fact that
$$
\|\phi_y\psi_y\rlap{--}\|\leq \|\phi_y\rlap{--}\|\|\psi_y\|_{L^{\infty}}
\leq C\|\phi_y\rlap{--}\|\|\psi_y\|^{\f12}\|\psi_{yy}\rlap{--}\|^{\f12}
\leq CN(\tau_1,\tau_2)(\|\psi_y\|+\|\psi_{yy}\rlap{--}\|).
$$
Substituting \eqref{psi_tau} into \eqref{jump-3} yields
$$
\ba
-\int_{\tau_1}^{\tau}[\psi_y\psi_{\tau}](\tau)d\tau
&\leq C\d\int_{\tau_1}^{\tau}\big(\d^{\f32}|\dot{\mb{X}}(\tau)|+\|\psi_{yy}\rlap{--}\|+\|(\phi_y,\psi_y,\z_y)\rlap{--}\|
+\d^{\f32}\big)^{\f12}\\
&\quad\times\|\psi_{\tau y}\rlap{--}\|^{\f12}
(\|\psi_{y}\|^{\f12}\|\psi_{yy}\rlap{--}\|^{\f12}+1)e^{-c(\tau-\tau_1)}d\tau\\
&\leq C\d\int_{\tau_1}^{\tau}\big(\|(\psi_{yy},\psi_{\tau y})\rlap{--}\|^2
+\|(\phi_y,\psi_y,\z_y)\rlap{--}\|^2+\d|\dot{\mb{X}}(\tau)|^2\big)d\tau+C\d.
\ea
$$
By Gagliardo-Nirenberg inequality and assumption \eqref{assumption}, it holds
$$
\ba
C\int_{\tau_1}^{\tau}\int_{\mathbb{R}}\hspace{-1.15em}-|\phi_y|\,|\psi_y|\,|\psi_{yy}|dyd\tau
&\leq C\int_{\tau_1}^{\tau}\|\phi_y\|\,\|\psi_y\|_{L^{\infty}}\|\psi_{yy}\rlap{--}\|d\tau\\
&\leq CN(\tau_1,\tau_2)\int_{\tau_1}^{\tau}\|\psi_y\|^{\f12} \|\psi_{yy}\rlap{--}\|^{\f32}d\tau\\
&\leq CN(\tau_1,\tau_2)\int_{\tau_1}^{\tau}(\|\psi_y\|^2+\|\psi_{yy}\rlap{--}\|^2)d\tau.
\ea
$$
Other terms can be treated similar as that of in Lemma \ref{le-1st}. Substituting the above estimates into \eqref{u-x} and using Lemmas \ref{le-L2}, \ref{le-1st},  we obtain
\be\label{es-u}
\ba
\|\psi_y(\tau)\|^2+\int_{\tau_1}^{\tau}\|\psi_{yy}\rlap{--}\|^2d\tau\leq C(N^2(\tau_1)+\d)+C\d\int_{\tau_1}^{\tau}\|\psi_{\tau y}\rlap{--}\|^2d\tau.
\ea
\ee


Next, we multiply \eqref{perturb}$_3$ by $-\z_{yy}$ to have
\be\label{t-1st}
\ba
&\f{R}{\gamma-1}\left(\f{\z_y^2}{2}\right)_\tau
-\f{R}{\gamma-1}(\z_y\z_\tau)_y+\f{\nu}{v}\z_{yy}^2\\
&=-\dot{\mathbf{X}}(\tau)\f{R\T^{-\mathbf{X}}_y}{\gamma-1}\z_{yy}
+p\psi_y\z_{yy}+(p-P^{-\mb{X}})U^{-\mb{X}}_y\z_{yy}
-\nu \T^{-\mb{X}}_{yy}\left(\f{1}{v}-\f{1}{V^{-\mb{X}}}\right)\z_{yy}\\
&+\f{\nu}{v^2}\phi_y\z_y\z_{yy}+\f{\nu}{v^2} \T^{-\mb{X}}_y\phi_y\z_{yy}
+\f{\nu}{v^2}V^{-\mb{X}}_y\z_y\z_{yy}+\nu V^{-\mb{X}}_y \T^{-\mb{X}}_y\left(\f{1}{v^2}-\f{1}{(V^{-\mb{X}})^2}\right)\z_{yy}\\
&-\f{1}{v}\psi_y^2\z_{yy}-\f{2}{v} U^{-\mb{X}}_y\psi_y\z_{yy}
-(U^{-\mb{X}}_y)^2\left(\f{1}{v}-\f{1}{V^{-\mb{X}}}\right)\z_{yy}.
\ea
\ee
Integrating \eqref{t-1st} with respect to $\tau$ and $y$ over $[\tau_1,\tau]\times\mathbb{R}_{\pm}$ gives
\be\label{t-x}
\ba
&\int_{\mathbb{R}}\f{R}{\gamma-1}\f{\z_y^2}{2}(\tau,y)dy
+\int_{\mathbb{R}}\hspace{-1.15em}-\f{\nu}{v}\f{\z_{yy}^2}{2}dy
\leq \int_{\mathbb{R}}\f{R}{\gamma-1}\f{\z_y^2}{2}(\tau_1,y)dy
-\f{R}{\gamma-1}\int_{\tau_1}^{\tau}[\z_y\z_{\tau}](\tau)d\tau\\
&+C\int_{\tau_1}^{\tau}|\dot{\mb{X}}(\tau)|^2\int_{\mathbb{R}}|\T^{-\mb{X}}_{y}|^2dyd\tau
+C\int_{\tau_1}^{\tau}\int_{\mathbb{R}}\psi_y^2dyd\tau
+C\int_{\tau_1}^{\tau}\int_{\mathbb{R}}\hspace{-1.15em}-|\phi_y|\,|\z_y|\,|\z_{yy}|dyd\tau\\
&+C\int_{\tau_1}^{\tau}\int_{\mathbb{R}}\hspace{-1.15em}-\psi_y^2|\z_{yy}|dyd\tau
+C\int_{\tau_1}^{\tau}\int_{\mathbb{R}}\hspace{-1.15em}-
|(U^{-\mb{X}}_{y},\T^{-\mb{X}}_{yy},V^{-\mb{X}}_y\T^{-\mb{X}}_y)|^2|(\phi,\z)|^2dyd\tau\\
&+C\int_{\tau_1}^{\tau}\int_{\mathbb{R}}\hspace{-1.15em}-|(V^{-\mb{X}}_y,U^{-\mb{X}}_{y},\T^{-\mb{X}}_y)|^2|(\phi_y,\psi_y,\z_y)|^2dyd\tau.
\ea
\ee
The jump across $y=0$ can be bounded as follows
\be\label{jump-4}
\ba
-\f{R}{\gamma-1}\int_{\tau_1}^{\tau}[\z_y\z_{\tau}](\tau)d\tau
&=-\f{R}{\gamma-1}\int_{\tau_1}^{\tau}\z_{\tau}(\tau,0)[\z_y](\tau)d\tau
=-\f{R}{\gamma-1}\int_{\tau_1}^{\tau}\z_{\tau}(\tau,0)[\t_y](\tau)d\tau\\
&=-\f{R}{\gamma-1}\int_{\tau_1}^{\tau}\z_{\tau}(\tau,0)\f{\t_y}{v}(\tau,0)[v](\tau)d\tau\\
&\leq C\int_{\tau_1}^{\tau}\|\z_{\tau}\|_{L^{\infty}}(\|\z_y\|_{L^{\infty}}+1)
[v](\tau_1)e^{-c(\tau-\tau_1)}d\tau\\
&\leq C\d\int_{\tau_1}^{\tau}\|\z_{\tau}\|^{\f12}\|\z_{\tau y}\rlap{--}\|^{\f12}
(\|\z_{y}\|^{\f12}\|\z_{yy}\rlap{--}\|^{\f12}+1)e^{-c(\tau-\tau_1)}d\tau.
\ea
\ee
By equation \eqref{perturb}$_3$ and assumption \eqref{assumption}, it holds
\be\label{z_tau}
\ba
\|\z_{\tau}\|&\leq C\big(|\dot{\mb{X}}(\tau)|\|\T^{-\mb{X}}_y\|+\|\z_{yy}\rlap{--}\|+\|\psi_y\|
+\|\phi_y\z_y\rlap{--}\|+\|\psi_y^2\|\\
&\quad+\|(U^{-\mb{X}}_y,(U^{-\mb{X}}_y)^2,V^{-\mb{X}}_y\T^{-\mb{X}}_y,\T^{-\mb{X}}_{yy})(\phi,\z)\rlap{--}\|
+\|(V^{-\mb{X}}_y,U^{-\mb{X}}_y,\T^{-\mb{X}}_y)(\phi_y,\psi_y,\z_y)\rlap{--}\|\big)
&\\
&\leq C\big(\d^{\f32}|\dot{\mb{X}}(\tau)|+\|(\psi_{yy},\z_{yy})\rlap{--}\|+\|(\phi_y,\psi_y,\z_y)\rlap{--}\|
+\d^{\f32}\big),
\ea
\ee
where we have used the facts that
$$
\|\phi_y\z_y\rlap{--}\|\leq \|\phi_y\rlap{--}\|\|\z_y\|_{L^{\infty}}
\leq C\|\phi_y\rlap{--}\|\|\z_y\|^{\f12}\|\z_{yy}\rlap{--}\|^{\f12}
\leq CN(\tau_1,\tau_2)(\|\z_y\|+\|\z_{yy}\rlap{--}\|)
$$
and
$$
\|\psi_y^2\|=\|\psi_y\|^2_{L^4}\leq C\|\psi_y\|^{\f32}\|\psi_{yy}\rlap{--}\|^{\f12}
\leq CN(\tau_1,\tau_2)\|\psi_y\|^{\f12}\|\psi_{yy}\rlap{--}\|^{\f12}
\leq CN(\tau_1,\tau_2)(\|\psi_y\|+\|\psi_{yy}\rlap{--}\|).
$$
Substituting \eqref{z_tau} into \eqref{jump-4} yields
$$
\ba
-\f{R}{\gamma-1}\int_{\tau_1}^{\tau}[\z_y\z_{\tau}](\tau)d\tau
&\leq C\d\int_{\tau_1}^{\tau}\big(\d^{\f32}|\dot{\mb{X}}(\tau)|+\|(\psi_{yy},\z_{yy})\rlap{--}\|+\|(\phi_y,\psi_y,\z_y)\rlap{--}\|
+\d^{\f32}\big)^{\f12}\\
&\quad\times\|\z_{\tau y}\rlap{--}\|^{\f12}
(\|\z_{y}\|^{\f12}\|\z_{yy}\rlap{--}\|^{\f12}+1)e^{-c(\tau-\tau_1)}d\tau\\
&\leq C\d\int_{\tau_1}^{\tau}\big(\|(\psi_{yy},\z_{yy},\z_{\tau y})\rlap{--}\|^2
+\|(\phi_y,\psi_y,\z_y)\rlap{--}\|^2+\d|\dot{\mb{X}}(\tau)|^2\big)d\tau+C\d.
\ea
$$
By Gagliardo-Nirenberg inequality and assumption \eqref{assumption}, it holds
$$
\ba
C\int_{\tau_1}^{\tau}\int_{\mathbb{R}}\hspace{-1.15em}-
|\phi_y|\,|\z_y|\,|\z_{yy}|dyd\tau
&\leq C\int_{\tau_1}^{\tau}\|\phi_y\rlap{--}\|\,\|\z_y\|_{L^{\infty}}\|\z_{yy}\rlap{--}\|d\tau\\
&\leq C\int_{\tau_1}^{\tau}N(\tau_1,\tau_2)\|\z_y\|^{\f12} \|\z_{yy}\rlap{--}\|^{\f32}d\tau\\
&\leq CN(\tau_1,\tau_2)\int_{\tau_1}^{\tau}(\|\z_y\|^2+\|\z_{yy}\rlap{--}\|^2)d\tau
\ea
$$
and
$$
\ba
C\int_{\tau_1}^{\tau}\int_{\mathbb{R}}\hspace{-1.15em}-\psi_y^2|\z_{yy}|dyd\tau
&\leq C\int_{\tau_1}^{\tau}\|\psi_y\|_{L^4}^2\|\z_{yy}\rlap{--}\|d\tau
\leq C\int_{\tau_1}^{\tau}\|\psi_y\|^{\f32} \|\psi_{yy}\rlap{--}\|^{\f12} \|\z_{yy}\rlap{--}\|d\tau\\
&\leq CN(\tau_1,\tau_2)\int_{\tau_1}^{\tau}\|\psi_y\|^{\f12} \|\psi_{yy}\rlap{--}\|^{\f12} \|\z_{yy}\rlap{--}\|d\tau\\
&\leq CN(\tau_1,\tau_2)\int_{\tau_1}^{\tau}(\|\psi_y\|^2+\|\psi_{yy}\rlap{--}\|^2+\|\z_{yy}\rlap{--}\|^2)d\tau.
\ea
$$
Other terms can be treated similar as that of in Lemma \ref{le-1st}. Substituting the above estimates into \eqref{t-x} and using Lemmas \ref{le-L2}, \ref{le-1st},  we have
\be\label{es-t}
\ba
\|\z_y(\tau)\|^2+\int_{\tau_1}^{\tau}\|\z_{yy}\rlap{--}\|^2d\tau
\leq C(N^2(\tau_1)+\d)+C\d\int_{\tau_1}^{\tau}\|(\z_{\tau y},\psi_{yy})\rlap{--}\|^2d\tau.
\ea
\ee
We combine \eqref{es-u} and \eqref{es-t} with suitably small $\d$ to have
\be\label{uy-1st}
\ba
\|(\psi_y,\z_y)(\tau)\|^2+\int_{\tau_1}^{\tau}\|(\psi_{yy},\z_{yy})\rlap{--}\|^2d\tau
\leq C(N^2(\tau_1)+\d)+C\d\int_{\tau_1}^{\tau}\|(\psi_{\tau y},\z_{\tau y})\rlap{--}\|^2d\tau.
\ea
\ee
Combination of \eqref{psi_tau}, \eqref{z_tau} and \eqref{uy-1st} gives
\be\label{L2-tau}
\ba
&\int_{\tau_1}^{\tau}\|(\psi_{\tau},\z_{\tau})\|^2d\tau\\
&\leq C\int_{\tau_1}^{\tau}\big(\d^3|\dot{\mb{X}}(\tau)|^2
+\|(\psi_{yy},\z_{yy})\rlap{--}\|^2+\|(\phi_y,\psi_y,\z_y)\rlap{--}\|^2
+G^s(\tau)+G_1(\tau)+G_2(\tau)\big)d\tau\\
&\leq C(N^2(\tau_1)+\d)+C\d\int_{\tau_1}^{\tau}\|(\psi_{\tau y},\z_{\tau y})\rlap{--}\|^2d\tau,
\ea
\ee
which together with \eqref{uy-1st} and suitably small $\d$ leads to the desired inequality \eqref{ut-1st} . The proof of Lemma \ref{le-1st-ut} is completed.

\hfill $\Box$

\begin{lemma}\label{le-2nd-yt}
Under the hypotheses of Proposition \ref{priori}, there exists constant $C>0$,
such that for any $\tau\in[\tau_1,\tau_2]$, it holds
\be\label{uzt}
\ba
\|(\psi_\tau,\z_\tau)(\tau)\|^2+\int_{\tau_1}^{\tau}\|(\psi_{\tau y},\z_{\tau y})\rlap{--}\hspace{0.07em}\|^2d\tau
\leq C(N^2(\tau_1)+\d).
\ea
\ee
\end{lemma}

\textbf{\emph{Proof}}: First, we apply the operator $\pa_\tau$ to the equation \eqref{perturb}$_2$ to have
$$
\ba
\psi_{\tau\tau}=\ddot{\mb{X}}(\tau)U^{-\mb{X}}_y-\big(s_3+\dot{\mb{X}}(\tau)\big)\dot{\mb{X}}(\tau)U^{-\mb{X}}_{yy}
+\left(\f{u_y}{v}-p\right)_{\tau y}-\Big(\f{U^{-\mb{X}}_y}{V^{-\mb{X}}}-P^{-\mb{X}}\Big)_{\tau y}.
\ea
$$
Multiplying the above equation by $\psi_\tau$, we have
\be\label{ut-1}
\ba
\left(\f{\psi_\tau^2}{2}\right)_{\tau}&+\f{\psi_{\tau y}^2}{v}
=\left(\psi_{\tau}\Big(\f{u_y}{v}-p\Big)_{\tau }-\psi_{\tau}\Big(\f{U^{-\mb{X}}_y}{V^{-\mb{X}}}-P^{-\mb{X}}\Big)_{\tau}
\right)_y\\
&+\ddot{\mb{X}}(\tau)\psi_{\tau }U^{-\mb{X}}_y-\big(s_3+\dot{\mb{X}}(\tau)\big)\dot{\mb{X}}(\tau)\psi_{\tau }U^{-\mb{X}}_{yy}\\
&+\psi_{\tau y}\f{u_y}{v^2}v_\tau-\psi_{\tau y}\f{U^{-\mb{X}}_{\tau y}}{v}
+\psi_{\tau y}\Big(\f{U^{-\mb{X}}_y}{V^{-\mb{X}}}\Big)_{\tau}
+\psi_{\tau y}(p-P^{-\mb{X}})_{\tau},
\ea
\ee
where
$$
\ba
&\f{u_y}{v^2}v_\tau-\f{U^{-\mb{X}}_{\tau y}}{v}+\Big(\f{U^{-\mb{X}}_y}{V^{-\mb{X}}}\Big)_{\tau}\\
=&\f{u_y^2}{v^2}-\f{(U^{-\mb{X}}_y)^2}{(V^{-\mb{X}})^2}
+\dot{\mb{X}}(\tau)\f{V^{-\mb{X}}_yU^{-\mb{X}}_y}{(V^{-\mb{X}})^2}
-U^{-\mb{X}}_{\tau y}\Big(\f{1}{v}-\f{1}{V^{-\mb{X}}}\Big)\\
=&\f{\psi_y^2}{v^2}+2\f{\psi_y}{v^2}U^{-\mb{X}}_y
+(U^{-\mb{X}}_y)^2\Big(\f{1}{v^2}-\f{1}{(V^{-\mb{X}})^2}\Big)
+\dot{\mb{X}}(\tau)\f{V^{-\mb{X}}_yU^{-\mb{X}}_y}{(V^{-\mb{X}})^2}
+\big(s_3+\dot{\mb{X}}(\tau)\big)U^{-\mb{X}}_{yy}\Big(\f{1}{v}-\f{1}{V^{-\mb{X}}}\Big)
\ea
$$
and
$$
\ba
&(p-P^{-\mb{X}})_{\tau}=\f{R\t_{\tau}}{v}-\f{R\T^{-\mb{X}}_{\tau}}{V^{-\mb{X}}}
-\f{pv_{\tau}}{v}+\f{P^{-\mb{X}}V^{-\mb{X}}_{\tau}}{V^{-\mb{X}}}\\
=&\f{R\z_{\tau}}{v}+R\T^{-\mb{X}}_{\tau}\Big(\f{1}{v}-\f{1}{V^{-\mb{X}}}\Big)
-\f{p u_y}{v}+\f{P^{-\mb{X}}}{V^{-\mb{X}}}(U^{-\mb{X}}_y-\dot{\mb{X}}(\tau)V^{-\mb{X}}_y)\\
=&\f{R\z_{\tau}}{v}-\f{p\psi_y}{v}-R\big(s_3+\dot{\mb{X}}(\tau)\big)\T^{-\mb{X}}_y\Big(\f{1}{v}-\f{1}{V^{-\mb{X}}}\Big)
-\dot{\mb{X}}(\tau)\f{P^{-\mb{X}}V^{-\mb{X}}_y}{V^{-\mb{X}}}
-\Big(\f{p}{v}-\f{P^{-\mb{X}}}{V^{-\mb{X}}}\Big)U^{-\mb{X}}_y.
\ea
$$
Integrating \eqref{ut-1} with respect to $\tau$ and $y$ over $[\tau_1,\tau]\times\mathbb{R}_{\pm}$, and using Cauchy inequality and the above two equalities, we obtain
\be\label{L2-ut}
\ba
&\int_{\mathbb{R}}\f{\psi_{\tau}^2}{2}(\tau,y)dy
+\int_{\tau_1}^{\tau}\int_{\mathbb{R}}\hspace{-1.15em}-
\f{\psi^2_{\tau y}}{2v}dyd\tau\\
&\leq\int_{\mathbb{R}}\f{\psi_{\tau}^2}{2}(\tau_1,y)dy
-\int_{\tau_1}^{\tau}\Big[\psi_{\tau}\Big(\f{u_y}{v}-p\Big)_{\tau }-\psi_{\tau}\Big(\f{U^{-\mb{X}}_y}{V^{-\mb{X}}}-P^{-\mb{X}}\Big)_{\tau}\Big](\tau)dy\\
&+\underbrace{\int_{\tau_1}^{\tau}\int_{\mathbb{R}}
|\ddot{\mb{X}}(\tau)|\,|\psi_{\tau }|\,|U^{-\mb{X}}_y|dyd\tau
+\int_{\tau_1}^{\tau}\int_{\mathbb{R}}(s_3+|\dot{\mb{X}}(\tau)|)
|\dot{\mb{X}}(\tau)|\,|\psi_{\tau }|\,|U^{-\mb{X}}_{yy}|dyd\tau}_{J_1}\\
&+\underbrace{C\int_{\tau_1}^{\tau}\int_{\mathbb{R}}\hspace{-1.15em}-
\Big(\psi_y^4+\psi_y^2|U^{-\mb{X}}_y|^2
+|U^{-\mb{X}}_y|^4\phi^2+|\dot{\mb{X}}(\tau)|^2|V^{-\mb{X}}_y U^{-\mb{X}}_y|^2
+(s_3+|\dot{\mb{X}}(\tau)|)^2|U^{-\mb{X}}_{yy}|^2\phi^2\Big)dyd\tau}_{J_2}\\ &+\underbrace{C\int_{\tau_1}^{\tau}\int_{\mathbb{R}}\hspace{-1.15em}-
\Big(\z_\tau^2+\psi_y^2+(s_3+|\dot{\mb{X}}(\tau)|)^2|\T^{-\mb{X}}_y|^2\phi^2
+|\dot{\mb{X}}(\tau)|^2|V^{-\mb{X}}_y|^2+|U^{-\mb{X}}_y|^2|(\phi,\z)|^2\Big)dyd\tau}_{J_3},
\ea
\ee
where the jump across $y=0$ in fact vanishes, i.e.,
$$
\ba
&\Big[\psi_{\tau}\Big(\f{u_y}{v}-p\Big)_{\tau }-\psi_{\tau}\Big(\f{U^{-\mb{X}}_y}{V^{-\mb{X}}}-P^{-\mb{X}}\Big)_{\tau}\Big](\tau)\\
=&[\psi_\tau](\tau)\Big(\f{u_y}{v}-p\Big)_{\tau}(\tau,0-)
+\psi_\tau(\tau,0+)\Big[\Big(\f{u_y}{v}-p\Big)_\tau\Big](\tau)
-[\psi_\tau](\tau)\Big(\f{U^{-\mb{X}}_y}{V^{-\mb{X}}}-P^{-\mb{X}}\Big)_{\tau}(\tau,0)\\
=&[\psi]_\tau(\tau)\Big(\f{u_y}{v}-p\Big)_{\tau}(\tau,0-)
+\psi_\tau(\tau,0+)\Big[\f{u_y}{v}-p\Big]_\tau(\tau)
-[\psi]_\tau(\tau)\Big(\f{U^{-\mb{X}}_y}{V^{-\mb{X}}}-P^{-\mb{X}}\Big)_{\tau}(\tau,0)=0.
\ea
$$
In the following, we control the remaining terms on the right-hand side of \eqref{L2-ut}. By \eqref{X}, Lemma \ref{le-shock} and \eqref{assumption}, it holds
$$
|\dot{\mathbf{X}}(\tau)|\leq \f{C}{\d}\|\psi(\tau)\|\|(V^{-\mathbf{X}}_y,U^{-\mathbf{X}}_y,\T^{-\mathbf{X}}_y)\|
\leq C\d^{\f12},
$$
and
$$
\ba
|\ddot{\mathbf{X}}(\tau)|&\leq \f{C}{\d}(s_3+|\dot{\mb{X}}(\tau)|)\int_{\mathbb{R}}|a^{-\mb{X}}_y|\,|(V^{-\mb{X}}_y,U^{-\mb{X}}_y,\T^{-\mb{X}}_y)|\,|\psi|dy\\
&\quad+\f{C}{\d}(s_3+|\dot{\mb{X}}(\tau)|)\int_{\mathbb{R}}|(V^{-\mb{X}}_{yy},U^{-\mb{X}}_{yy},\T^{-\mb{X}}_{yy},(V^{-\mb{X}}_{y})^2,(\T^{-\mb{X}}_{y})^2)|\,|\psi|dy\\
&\quad+\f{C}{\d}\int_{\mathbb{R}}|(V^{-\mb{X}}_y,U^{-\mb{X}}_y,\T^{-\mb{X}}_y)|\,|\psi_{\tau}|dy\\
&\leq C\d^{\f12}\big(\sqrt{G^s(\tau)}+\|\psi_{\tau}(\tau)\|\big).
\ea
$$
Then, we use the above two inequalities, \eqref{basic} and \eqref{L2-tau} to have
$$
\ba
J_1&\leq C\d^{\f12}\int_{\tau_1}^{\tau}\big(\sqrt{G^s(\tau)}+\|\psi_{\tau}\|\big)
\|\psi_{\tau}\|\,\|U^{-\mb{X}}_y\|d\tau
+C\int_{\tau_1}^{\tau}|\dot{\mb{X}}(\tau)|\|\psi_{\tau}\|\,\|U^{-\mb{X}}_{yy}\|d\tau\\
&\leq C\d^2\int_{\tau_1}^{\tau}\big(G^s(\tau)+\|\psi_\tau\|^2+\d|\dot{\mb{X}}(\tau)|^2\big)d\tau
\leq C(N^2(\tau_1)+\d)+C\d\int_{\tau_1}^{\tau}\|(\psi_{\tau y},\z_{\tau y})\rlap{--}\|^2d\tau.
\ea
$$
Similarly, using Sobolev inequality, \eqref{basic} and \eqref{ut-1st}, we have
$$
\ba
J_2+J_3&\leq C\int_{\tau_1}^{\tau}\big(\|\psi_{yy}\rlap{--}\|^2+\|\z_{\tau}\|^2
+\|\psi_y\|^2+G^s(\tau)+G_1(\tau)+G_2(\tau)+\d|\dot{\mb{X}}(\tau)|^2\big)d\tau\\
&\leq C(N^2(\tau_1)+\d)+C\d\int_{\tau_1}^{\tau}\|(\psi_{\tau y},\z_{\tau y})\rlap{--}\|^2d\tau.
\ea
$$
Substituting the above estimates into \eqref{L2-ut} with suitably small $\d$ gives
\be\label{es-ut}
\|\psi_{\tau}(\tau)\|^2+\int_{\tau_1}^{\tau}\|\psi_{\tau y}\rlap{--}\|^2d\tau
\leq C(N^2(\tau_1)+\d)+C\d\int_{\tau_1}^{\tau}\|\z_{\tau y}\rlap{--}\|^2d\tau.
\ee
Next, we apply the operator $\pa_\tau$ to the equation \eqref{perturb}$_3$ to have
$$
\ba
\f{R\z_{\tau\tau}}{\gamma-1}&=\ddot{\mb{X}}(\tau)\f{R\T^{-\mb{X}}_y}{\gamma-1}
-\big(s_3+\dot{\mb{X}}(\tau)\big)\dot{\mb{X}}(\tau)\f{R\T^{-\mb{X}}_{yy}}{\gamma-1}
+\nu\Big(\f{\t_y}{v}\Big)_{\tau y}-\nu\Big(\f{\T^{-\mb{X}}_y}{V^{-\mb{X}}}\Big)_{\tau y}\\
&+\left(u_y\Big(\f{u_y}{v}-p\Big)\right)_\tau
-\Big(U^{-\mb{X}}_y\Big(\f{U^{-\mb{X}}_y}{V^{-\mb{X}}}-P^{-\mb{X}}\Big)\Big)_\tau.
\ea
$$
Multiplying the above equation by $\z_\tau$, we have
$$
\ba
\f{R}{\gamma-1}\Big(\f{\z_\tau^2}{2}\Big)^2&+\nu\f{\z_{\tau y}^2}{v}
=\nu\left(\z_\tau\Big(\f{\t_y}{v}\Big)_\tau-\z_\tau\Big(\f{\T^{-\mb{X}}_y}{V^{-\mb{X}}}\Big)_{\tau}\right)_y\\
&+\ddot{\mb{X}}(\tau)\z_\tau\f{R\T^{-\mb{X}}_y}{\gamma-1}
-\big(s_3+\dot{\mb{X}}(\tau)\big)\dot{\mb{X}}(\tau)\z_{\tau}\f{R\T^{-\mb{X}}_{yy}}{\gamma-1}
-\nu\z_{\tau y}\f{\T^{-\mb{X}}_{\tau y}}{v}+\nu\z_{\tau y}\f{\t_y}{v^2}v_\tau\\
&+\nu\z_{\tau y}\Big(\f{\T^{-\mb{X}}_y}{V^{-\mb{X}}}\Big)_\tau
+\z_\tau u_{\tau y}\Big(\f{u_y}{v}-p\Big)+\z_\tau u_y\Big(\f{u_y}{v}-p\Big)_\tau\\
&-\z_\tau U^{-\mb{X}}_{\tau y}\Big(\f{U^{-\mb{X}}_y}{V^{-\mb{X}}}-P^{-\mb{X}}\Big)
-\z_\tau U^{-\mb{X}}_{ y}\Big(\f{U^{-\mb{X}}_y}{V^{-\mb{X}}}-P^{-\mb{X}}\Big)_\tau.
\ea
$$
Integrating the above equality with respect to $\tau$ and $y$ over $[\tau_1,\tau]\times\mathbb{R}_{\pm}$ gives
\be\label{zt}
\ba
&\int_{\mathbb{R}}\f{R}{\gamma-1}\f{\z_\tau^2}{2}(\tau,y)dy
+\int_{\tau_1}^{\tau}\int_{\mathbb{R}}\hspace{-1.15em}-\nu\f{\z_{\tau y}^2}{v}dyd\tau\\
&=\int_{\mathbb{R}}\f{R}{\gamma-1}\f{\z_\tau^2}{2}(\tau_1,y)dy
-\int_{\tau_1}^{\tau}\nu\Big[\z_\tau\Big(\f{\t_y}{v}\Big)_\tau-\z_\tau\Big(\f{\T^{-\mb{X}}_y}{V^{-\mb{X}}}\Big)_\tau\Big](\tau)d\tau\\
&+\int_{\tau_1}^{\tau}\int_{\mathbb{R}}\Big(
\ddot{\mb{X}}(\tau)\z_\tau\f{R\T^{-\mb{X}}_y}{\gamma-1}
-\big(s_3+\dot{\mb{X}}(\tau)\big)\dot{\mb{X}}(\tau)\z_{\tau}\f{R\T^{-\mb{X}}_{yy}}{\gamma-1}
\Big)dyd\tau\\
&+\int_{\tau_1}^{\tau}\int_{\mathbb{R}}\hspace{-1.15em}-\Big\{
-\nu\z_{\tau y}\f{\T^{-\mb{X}}_{\tau y}}{v}+\nu\z_{\tau y}\f{\t_y}{v^2}v_\tau
+\nu\z_{\tau y}\Big(\f{\T^{-\mb{X}}_y}{V^{-\mb{X}}}\Big)_\tau
+\z_\tau u_{\tau y}\Big(\f{u_y}{v}-p\Big)\\
&+\z_\tau u_y\Big(\f{u_y}{v}-p\Big)_\tau-\z_\tau U^{-\mb{X}}_{\tau y}\Big(\f{U^{-\mb{X}}_y}{V^{-\mb{X}}}-P^{-\mb{X}}\Big)-\z_\tau U^{-\mb{X}}_{ y}\Big(\f{U^{-\mb{X}}_y}{V^{-\mb{X}}}-P^{-\mb{X}}\Big)_\tau
\Big\}dyd\tau,
\ea
\ee
where the jump in fact vanishes, i.e.,
\be
\ba
&\Big[\z_\tau\Big(\f{\t_y}{v}\Big)_\tau-\z_{\tau}\Big(\f{\T^{-\mb{X}}_y}{V^{-\mb{X}}}\Big)_\tau\Big](\tau)\\
=&[\z_\tau](\tau)\Big(\f{\t_y}{v}\Big)_{\tau}(\tau,0-)+\z_\tau(\tau,0+)\Big[\Big(\f{\t_y}{v}\Big)_{\tau}\Big](\tau)
-[\z_\tau](\tau)\Big(\f{\T^{-\mb{X}}_y}{V^{-\mb{X}}}\Big)_{\tau}(\tau,0)\\
=&[\z]_{\tau}(\tau)\Big(\f{\t_y}{v}\Big)_{\tau}(\tau,0-)+\z_\tau(\tau,0+)\Big[\f{\t_y}{v}\Big]_{\tau}(\tau)
-[\z]_{\tau}(\tau)\Big(\f{\T^{-\mb{X}}_y}{V^{-\mb{X}}}\Big)_{\tau}(\tau,0)=0.
\ea
\ee
The other terms on the right-hand side of \eqref{zt} can be controlled as \eqref{es-ut}, therefore, it holds
\be\label{es-zt}
\|\z_{\tau}(\tau)\|^2+\int_{\tau_1}^{\tau}\|\z_{\tau y}\rlap{--}\|^2d\tau
\leq C(N^2(\tau_1)+\d)+C\d\int_{\tau_1}^{\tau}\|\psi_{\tau y}\rlap{--}\|^2d\tau.
\ee
Thus, combination of \eqref{es-ut} and \eqref{es-zt} together with suitably small $\d$ gives the desired inequality \eqref{uzt}. The proof of Lemma \ref{le-2nd-yt} is completed.

\hfill $\Box$

Finally, combination of Lemmas \ref{le-L2}-\ref{le-2nd-yt} yields Proposition \ref{priori}.

\section*{Acknowledgments}
The authors would like to thank the anonymous referee for his/her helpful comments and suggestions, which greatly improve the representation of the manuscript.
The research of F. M. Huang is partially supported by National Key R\&D Program of China (No. 2021YFA1000800), the NSFC of China (Nos. 11371349 and 11688101).
The research of T. Wang is partially supported by the National Key R\&D Program of China (No. 2022YFA1007700), the NSFC of China (Nos. 11971044 and 12371215).


\small
\begin{thebibliography}{99}

\bibitem{BB}
S. Bianchini and A. Bressan, \emph{Vanishing viscosity solutions of nonlinear hyperbolic systems}. Ann. of Math. \textbf{161} (2005), 223-342.

\bibitem{B}
A. Bressan, \emph{BV-solutions to hyperbolic systems by vanishing viscosity}. S.I.S.S.A., Trieste, Italy, 2000.

\bibitem{BY}
A. Bressan and T. Yang, \emph{On the convergence rate of vanishing viscosity approximations}. Comm. Pure Appl. Math. \textbf{57} (2004), 1075-1109.

\bibitem{BHWY}
A. Bressan, F. M. Huang, Y. Wang and T. Yang,
\emph{On the convergence rate of vanishing viscosity approximations for nonlinear hyperbolic systems}.
SIAM J. Math. Anal. \textbf{44} (2012), 3537-3563.


\bibitem{Chen-acta}
G.-Q. Chen, \emph{Convergence of the Lax-Friedrichs scheme for isentropic gas dynamics (III)}.
Acta Math. Sci. \textbf{6B} (1986), 75-120 (in English); \textbf{8A} (1988), 243-276 (in Chinese).

\bibitem{Chen-proc}
G.-Q. Chen, \emph{Remarks on DiPerna's paper: ``Convergence of the viscosity method for isentropic gas dynamics" [Comm. Math. Phys. 91 (1983), 1-30]}.
Proc. Amer. Math. Soc. \textbf{125} (1997), 2981-2986.

\bibitem{CP}
G. Chen and M. Perepelitsa, \emph{Vanishing viscosity limit of the Navier-Stokes equations to the Euler equations for compressible fluid flow}. Comm. Pure Appl. Math. \textbf{63} (2010),  1469-1504.




\bibitem{DCL}
X. Ding, G.-Q. Chen and P. Luo, \emph{Convergence of the Lax-Friedrichs scheme for the isentropic gas dynamics (I)-(II)}.
Acta Math. Sci. \textbf{5B}, 483-500 (1985), 501-540 (in English);
\textbf{7A} (1987), 467-480;
\textbf{8A} (1989), 61-94 (in Chinese);
\emph{Convergence of the fractional step Lax-Friedrichs scheme and Godunov scheme for the isentropic system of gas dynamics}.
Comm. Math. Phys. \textbf{121} (1989), 63-84.


\bibitem{Diperna}
R. J. DiPerna, \emph{Convergence of the viscosity method for isentropic gas dynamics}. Comm. Math. Phys.  \textbf{91} (1983), 1-30.

\bibitem{Gilbarg}
D. Gilbarg, \emph{The existence and limit behavior of the one-dimensional shock layer}. Amer. J. Math. \textbf{73} (1951), 256-274.

\bibitem{Goodman-Xin}
J. Goodman and Z. P. Xin, \emph{Viscous limits for piecewise smooth solutions to systems of conservation laws}. Arch. Ration. Mech. Anal. \textbf{121} (1992), 235-265.


\bibitem{Hoff}
D. Hoff, \emph{Discontinuous solutions of the Navier-Stokes equations for
compressible flow}. Arch. Ration. Mech. Anal. \textbf{114} (1991), 15-46.


\bibitem{HL-1}
D. Hoff and T. P. Liu, \emph{The inviscid limit for the Navier-Stokes equations of compressible, isentropic flow with shock data}.
Indiana Univ. Math. J \textbf{38} (1989), 861-915.



%
%
%
\bibitem{Huang-matsumura}
 F. M. Huang and A. Matsumura,
\emph{ Stability of a composite wave of two viscous shock waves for the full compressible Navier-Stokes equation},
Comm. Math. Phys. \textbf{289} (2009), 841-861.
%
%
%
%
%

\bibitem{HJW}
F. M. Huang, S. Jiang and Y. Wang,
\emph{Zero dissipation limit of full compressible Navier-Stokes equations with a Riemann initial data}.
Commun. Inf. Syst. \textbf{13} (2013), 211-246.

\bibitem{HWY-KRM}
F. M. Huang, Y. Wang and T. Yang, \emph{Fluid dynamic limit to the Riemann solutions of Euler equations: I. superposition of rarefaction waves and contact discontinuity}.  Kinet. Relat. Models \textbf{3} (2010), 685-728.



\bibitem{HWY-ARMA}
F. M. Huang, Y. Wang and T. Yang, \emph{Vanishing viscosity limit of the compressible Navier-Stokes equations for solutions to Riemann problem}.
Arch. Ration. Mech. Anal. \textbf{203} (2012), 379-413.

\bibitem{HWWY}
F. M. Huang, Y. Wang, Y. Wang, T. Yang, \emph{The limit of the Boltzmann equation to the Euler equations for Riemann problems}.
SIAM J. Math. Anal. \textbf{45} (2013), 1741-1811.

\bibitem{Huang-WangZ}
F. M. Huang and Z. Wang, \emph{Convergence of viscosity solutions for isothermal gas dynamics}.
SIAM J. Math. Anal. \textbf{34} (2002), 595-610.




\bibitem{JNS}
S. Jiang, G. X. Ni and W. J. Sun, \emph{Vanishing viscosity limit to rarefaction waves for the Navier-Stokes equations of one-dimensional compressible heat-conducting fluids}. SIAM J. Math. Anal. \textbf{38} (2006), 368-384.



\bibitem{KV-JEMS}
M.-J. Kang and A. Vasseur, \emph{Contraction property for large perturbations of shocks of the barotropic Navier-Stokes system}.
J. Eur. Math. Soc.  \textbf{23} (2021) 585-638.



\bibitem{KV-Invention}
M.-J. Kang and A. Vasseur, \emph{Uniqueness and stability of entropy shocks to the isentropic Euler system in a class of inviscid limits from a large family of Navier-Stokes systems}.
Invent. math.  \textbf{224} (2021),  55-146.



\bibitem{KVW}
M.-J. Kang,  A. Vasseur and Y. Wang, \emph{Time-asymptotic stability of composite waves of viscous shock wave and rarefaction for barotropic Navier-Stokes equations},  Adv. Math. \textbf{419} (2023), Paper No. 108963, 66 pp.

\bibitem{KVW-2}
M.-J. Kang,  A. Vasseur and Y. Wang, \emph{Time-asymptotic stability of generic Riemann solutions for compressible Navier-Stokes-Fourier equations}, https://arxiv.org/abs/2306.05604.

\bibitem{KM}
S. Kawashima and A. Matsumura, \emph{Asymptotic stability of traveling wave solutions of systems for one-dimensional gas motion},
Comm. Math. Phys. \textbf{101} (1985), 97-127.


\bibitem{Lax}
P. Lax, \emph{Hyperbolic systems of conservation laws, II}. Comm. Pure Appl. Math. {\bf 10} (1957), 537-566.

\bibitem{Lions-PS}
P.-L. Lions, B. Perthame and P. E. Souganidis, \emph{Existence and stability of entropy solutions for the hyperbolic systems of isentropic gas dynamics in Eulerian and Lagrangian coordinates}.
Comm. Pure Appl. Math. \textbf{49} (1996), 599-638.

\bibitem{Lions-PT}
P.-L.Lions, B. Perthame and E. Tadmor,
\emph{Kinetic formulation of the isentropic gas dynamics and p-systems}.
Comm. Math. Phys. \textbf{163} (1994), 415-431.

\bibitem{Ma}
S. X. Ma, \emph{Zero dissipation limit to strong contact discontinuity for the 1-D compressible Navier-Stokes equations}.
J Differential Equations \textbf{258} (2010), 95-110.

%


%
%
%
%
%
%
%
%
%

\bibitem{MN-92}
A. Matsumura, K. Nishihara, \emph{Global stability of the rarefaction wave of a one-dimensional model system for compressible viscous gas}, Comm. Math. Phys. \textbf{144} (1992), 325-335.

\bibitem{SS}
M. R. I. Schrecker, S. Schulz, \emph{Vanishing viscosity limit of the compressible Navier-Stokes equations with general pressure law}.
SIAM J. Math. Anal. \textbf{51} (2019), 2168-2205.

\bibitem{stokes}
G. G. Stokes, \emph{On a difficulty in the theory of sound}. Philos. Mag. \textbf{33} (1848), 349-356.



\bibitem{WW}
T. Wang and Y. Wang,
\emph{Nonlinear stability of planar viscous shock wave to three-dimensional compressible Navier-Stokes equations}, to appear in  J. Eur. Math. Soc., https://arxiv.org/abs/2204.09428.



%
%

\bibitem{WY}
Y. Wang, \emph{Zero dissipation limit of the compressible heat-conducting Navier-Stokes equations in the presence of the shock}.
Acta Mathematica Scientia \textbf{28} (2008), 727-748.

\bibitem{Xin}
Z. P. Xin, \emph{Zero dissipation limit to rarefaction waves for the one-dimensional Navier-Stokes equations of compressible isentropic gases}.
Comm. Pure Appl. Math. \textbf{46} (1993), 621-665.
%

\bibitem{Xin-Z}
Z. P. Xin and H. H. Zeng, \emph{Convergence to the rarefaction waves for the nonlinear Boltzmann equation and compressible Navier-Stokes equations}.
J. Differential Equations. \textbf{249} (2010), 827-871.

\bibitem{Yu-ARMA}
S. H. Yu, \emph{Zero-dissipation limit of solutions with shocks for systems of hyperbolic conservation laws}. Arch. Ration. Mech. Anal. \textbf{146} (1999), 275-370.



\bibitem{ZPT}
Y. H. Zhang, R. H. Pan and Z. Tan, \emph{Zero dissipation limit to a Riemann solution consisting of two shock waves for the 1D compressible isentropic Navier-Stokes equations}. Sci. China Math. \textbf{56} (2013), 2205-2232.


\bibitem{ZPWT}
Y. H. Zhang, R. H. Pan, Y. Wang and Z. Tan,
\emph{Zero dissipation limit with two interacting shocks of the 1D non-isentropic Navier-Stokes equations}.
Indiana Univ. Math. J. \textbf{62} (2013), 249-309.




\end{thebibliography}
\end{document}